\documentclass[12pt,oneside,british,a4wide]{amsart}
\usepackage{amsmath, amssymb,verbatim,appendix}
\usepackage[pagebackref=true]{hyperref}
\usepackage[mathscr]{eucal}
\usepackage{amscd}
\usepackage{amsthm}
\usepackage{stmaryrd}
\usepackage{enumitem}
\usepackage{comment}
\usepackage[latin1]{inputenc} 
\usepackage{tikz}
\usetikzlibrary{shapes,arrows}
\usepackage{cite}
\usepackage{url}
\usepackage{float}

\usepackage{setspace,a4wide,color,xcolor,graphicx}

\def\showauthornotes{1}

\ifnum\showauthornotes=1
\newcommand{\Authornote}[2]{{\sf\small\color{red}{[#1: #2]}}}
\else
\newcommand{\Authornote}[2]{}
\fi

\newtheorem{theorem}{Theorem}[section]
\newtheorem{lemma}[theorem]{Lemma}
\newtheorem{corollary}[theorem]{Corollary}
\newtheorem{proposition}[theorem]{Proposition}
\newtheorem{claim}[theorem]{Claim}
\newtheorem{fact}[theorem]{Fact}

\numberwithin{equation}{section}
\setlist[enumerate]{font={\rmfamily}}
\setlist[enumerate,1]{label={(\roman*)}}

\theoremstyle{definition}

\newtheorem{definition}[theorem]{Definition}
\newtheorem{remark}[theorem]{Remark}

\def\F{\mathbb{F}}

\def\FOP{\mathrm{FOP}}
\def\NFOP{\mathrm{NFOP}}
\def\IP{\mathrm{IP}}
\def\NIP{\mathrm{NIP}}

\def\VC{\mathrm{VC}}

\def\e{\epsilon}

\newcommand{\bbar}{\bar{b}}

\newcommand{\calL}{\mathcal{L}}
\newcommand{\calF}{\mathcal{F}}

\newcommand{\calB}{\mathcal{B}}

\newcommand{\calP}{\mathcal{P}}

\newcommand{\calQ}{\mathcal{Q}}
\newcommand{\calR}{\mathcal{R}}
\newcommand{\calS}{\mathcal{S}}

\newcommand{\calX}{\mathcal{X}}

\newcommand{\bY}{\mathbf{Y}}
\newcommand{\bH}{\mathbf{H}}
\newcommand{\bR}{\mathbf{R}}

\newcommand{\codim}{{\mathrm{codim}}}

\newcommand{\Red}{{\mathrm{Red}}}
\newcommand{\Span}{{\mathrm{Span}}}
\newcommand{\At}{{\mathrm{At}}}
\newcommand{\rk}{{\mathrm{rk}}}
\newcommand{\Stab}{{\mathrm{Stab}}}

\def\Red{\operatorname{Red}}

\newenvironment{proofof}[1]{\indent{\itshape Proof of #1}.\;}{\qed}

\begin{document}

\title[On the quadratic complexity of subsets of $\F_p^n$ of bounded $\mathrm{VC_{2}}$-dimension]{On the quadratic complexity of \\subsets of $\F_p^n$ of bounded $\mathrm{VC_{2}}$-dimension}

\author{C. Terry}

\author{J. Wolf}

\address{Department of Mathematics, Statistics, and Computer Science, University of Illinois Chicago, Chicago IL 60607, USA}

\email{caterry@uic.edu}

\address{Department of Pure Mathematics and Mathematical Statistics, Centre for Mathematical Sciences, Wilberforce Road, Cambridge CB3 0WB, UK}

\email{julia.wolf@dpmms.cam.ac.uk}

\date{}

\begin{abstract}
In prior work, we showed that subsets of $\F_{p}^{n}$ of $\mathrm{VC_{2}}$-dimension at most $k$ are well approximated by a union of atoms of a quadratic factor of complexity $(\ell,q)$, where the complexity $\ell$ of the linear part and the complexity $q$ of the quadratic part are both bounded in terms of $k$, $p$, and the desired level of approximation $\mu$. A key tool in the proof of this result was an arithmetic regularity lemma for the Gowers $U^3$-norm by Green and Tao, which resulted in tower-type bounds (in terms of $\mu^{-1}$) on both $\ell$ and $q$. 

In the present paper we show that for sets of bounded $\VC_2$-dimension, the bound on $q$ can be substantially improved. Specifically, we will prove that any set $A\subseteq G=\F_p^n$ of $\VC_2$-dimension at most $k$ is approximately equal (up to error $\mu |G|$) to a union of atoms of a quadratic factor whose quadratic complexity is at most $\log_p(\mu^{-k-o(1)})$,  implying that the purely quadratic component of the factor partitions the group into $\mu^{-k-o(1)}$ many parts. 

We achieve this by using our earlier result to obtain an initial quadratic factor $\calB$, and then applying a generalization of an argument of Alon, Fox and Zhao for subsets of $\F_{p}^{n}$ of bounded $\mathrm{VC}$-dimension to the label space (also known as ``configuration space") of $\calB$.  A related strategy was employed in earlier work of the authors on $\NFOP_2$ subsets of $\F_p^n$, and in work of the first author in the context of 3-uniform hypergraphs.
\end{abstract}

\maketitle

\tableofcontents

\section{Introduction}

In part as a result of a fruitful interplay with ideas from model theory, there has been renewed interest in recent years in regularity lemmas as results in their own right, as opposed to mere tools for obtaining the kinds of applications for which they were originally conceived. This applies in particular to regularity lemmas for subsets of finite abelian groups, known as \emph{arithmetic regularity lemmas}, whose study was initiated by Green in \cite{Green.2005}.

In the context of the elementary abelian $p$-group $\F_p^n$, which continues to enjoy popularity as a toy setting in additive combinatorics, such a regularity lemma roughly says that given a subset $A\subseteq \F_p^n$, there exists a subspace $H\leqslant \F_p^n$ of bounded codimension with the property that $A$ is ``uniform" on almost all cosets of $H$. A notion of uniformity that is useful for (at least certain) applications can be defined in terms of the Fourier transform of the indicator function of $A$ restricted to cosets of $H$, or (equivalently) in terms of the \emph{local $U^2$ norms} defined in \cite{Terry.2021d}. In modern parlance, the regularity partition of $\F_p^n$ into cosets of $H$ is often referred to as a \emph{linear factor}, whose complexity can be taken to be the codimension of $H$ (see Section \ref{subsec:linquadfac} for a recap). For the purposes of this introduction,  it suffices to say that under this identification, the cosets of $H$ are instead referred to as the \emph{atoms of $\calL$}, each of which is written as $L(a)$ for some $a\in \F_p^{\ell}$, where $\ell=\codim(H)$.

With this terminology, a more formal statement of Green's \emph{linear} arithmetic regularity lemmas is Theorem \ref{thm:linearreg} below.

\begin{theorem}[$U^2$ arithmetic regularity \cite{Green.2005}]\label{thm:linearreg}
For every prime $p$ and $\e\in (0,1)$, there exists $C=C(\e)$ such that the following holds for all sufficiently large $n$.
For any $A\subseteq \F_p^n$, there is a linear factor $\calL$ of complexity $\ell\leq C$ and a set $\Gamma\subseteq \F_p^{\ell}\times \F_p^{\ell}$ such that
\begin{enumerate}[label=\normalfont(\roman*)] 
\item $|\Gamma|\geq (1-\e) p^{2\ell}$;
\item for all $d=(d_1,d_2)\in \Gamma$, $\|1_A-\alpha_{\Sigma(d)}\|_{U^2(d)}<\e$, where $\alpha_{\Sigma(d)}$ is the density of $A$ on the atom $L(\Sigma(d))=L(d_1+d_2)$.
\end{enumerate}
\end{theorem}

For the narrative of this introduction, the precise definitions of the terms used in Theorem \ref{thm:linearreg} are not strictly necessary. It suffices to say that, in the notation of Theorem \ref{thm:linearreg}, the atoms of the form $L(\Sigma(d))$ in point (ii) are considered \emph{uniform with respect to $A$}, and the lower bound on $\Gamma$ in (i) ensures almost all atoms of $\calL$ are uniform in this sense.

The proof in \cite{Green.2005} produced a bound on $C(\e)$ which is (in substance) a tower of height $\e^{-3}$. In the same paper, Green also showed that a bound of this shape is, in general, necessary (his lower bound was subsequently improved in \cite{Hosseini.2016} to one that is essentially tight).

Theorem \ref{thm:linearreg} was modelled on a celebrated regularity lemma for graphs of Szemer\'edi \cite{Szemeredi.1975}, which provides a way of partitioning large graphs into a bounded number of components, almost all of which are ``uniform" (suitably defined). Graph regularity enjoys similarly poor quantitative features, as was famously proved by Gowers \cite{Gowers.1997}.  However, it was shown about 20 years ago \cite{Alon.2007,Lovasz.2007} that under certain hypotheses, the upper bound on the size of the regularity partition can be substantially improved. More recently, this phenomenon was also established in the context of subsets of groups. The first such result was an arithmetic regularity lemma for subsets of $\F_p^n$ under the assumption of \emph{stability} \cite{Terry.2019}, which was followed by work of several authors in the more general setting of \emph{bounded $\VC$-dimension}. We recall the notion of $\VC$-dimension as applied to a subset of an abelian group, which will be the focus of this paper. 

\begin{definition}[$\VC$-dimension of a subset of a group]\label{def:vcdim}
Given an integer $k\geq 1$, an abelian group $G$ and a subset $A\subseteq G$, we say that $A$ \emph{has $\VC$-dimension at least $k$} (or \emph{has $k$-$\IP$})\footnote{$\IP$ stands for the ``independence property'' in model theory.} if there exist elements $\{a_i: i\in [k]\}\cup \{b_S:S\subseteq [k]\}$ in $G$ such that $a_i+ b_S\in A$ if and only if $i\in S$.  

The \emph{$\VC$-dimension of $A$}, denoted $\VC(A)$, is defined to be the largest integer $k\geq 0$ such that $A$ has $\VC$-dimension at least $k$, when such a $k$ exists, and is defined to be $0$ otherwise.
\end{definition}
 
A set $A\subseteq G$ has $\VC$-dimension at least $k$ if and only if the set system $\{A-x:x\in G\}$ on $G$ shatters a set of size $k$, meaning there exists $T\subseteq G$ of size $k$, and for any subset $S\subseteq T$, an element $b_{S}\in G$ such that $(A-b_{S})\cap T=S$.   It is easy to see that the only subsets of $G$ with VC-dimension $0$ are the empty set and the whole group $G$.  

Generalising \cite{Terry.2019}, Alon, Fox and Zhao \cite{Alon.2018is} proved the following structure theorem for subsets of $\F_p^n$ of bounded $\VC$-dimension.

\begin{theorem}[Structure theorem for sets of bounded $\VC$-dimension \cite{Alon.2018is}]\label{thm:vc}
Let $p$ be a prime, let $k\geq 0$ be an integer, and let $\e\in (0,1)$.  For any $A\subseteq \F_p^n$ satisfying  $\VC(A)\leq k$, there exist  a subgroup $H\leqslant \F_p^n$ of index at most $\e^{-k-o_{k,p}(1)}$ and a union $Y$ of cosets of $H$ such that $|A\Delta Y|\leq \e|G|$, where $o_k(1)$ tends to $0$ as $\e$ tends to $0$, at a rate depending on $k$ and $p$. 
\end{theorem}
 
For results of a similar flavour (but in a more general context), see \cite{Sisask.2018,Conant.2018zd, Conant.2025}.  It is an immediate corollary of Theorem \ref{thm:vc} that the density of the set $A$ is near $0$ or $1$ on almost all cosets of the set $H$.  Restating this in terms of linear factors yields the following corollary (recall that a subgroup of index $m$ corresponds to a linear factor of complexity $\log_p(m)$).

\begin{corollary}[$U^2$ arithmetic regularity for sets of bounded VC-dimension]\label{cor:linearregvc}
For all primes $p>2$ and integers $k\geq 0$, there exists $C=C(k)$ such that the following holds for all $\e\in (0,1)$ and all sufficiently large $n$.

Given $A\subseteq \F_p^n$ with $\VC(A)\leq k$, there is a linear factor $\calL$ of complexity $\ell$ and a set $\Gamma\subseteq \F_p^{\ell}\times \F_p^{\ell}$ such that
\begin{enumerate}[label=\normalfont(\roman*)]
\item the complexity of $\calL$ satisfies $\ell \leq \log_p(\e^{-C})$;
\item $|\Gamma|\geq (1- \e) p^{2\ell}$;
\item for all $d\in \Gamma$, $\alpha_{\Sigma(d)}\in [0,\e)\cup (1-\e,1]$.

\end{enumerate}
\end{corollary}

Moreover, it follows from (a version of) \cite[Proposition 1]{Terry.2019} and \cite[Lemma 5]{Terry.2019} that in the context of sets of bounded $\VC$-dimension, having density near $0$ or $1$ is approximately equivalent to being ``uniform" in the sense of Theorem \ref{thm:linearreg} . 
 
\begin{fact}[Local homogeneity and local $U^2$ uniformity are approximately equivalent \cite{Terry.2019}]\label{fct:unihom}
Let $p$ be a fixed prime. For all $\e>0$ and all integers $k\geq 1$, $\ell\geq 1$, and $n > \ell$ the following holds.  Let $\calL$ be a linear factor on $\F_p^n$ of complexity $\ell$, let $d=(d_1,d_2)\in \F_p^\ell\times \F_p^{\ell}$ and suppose that $A\subseteq \F_p^n$ has density $\alpha_{\Sigma(d)}$ on the atom $L(\Sigma(d))=L(d_1+d_2)$. 

Then $\alpha_{\Sigma(d)}\in [0,\e)\cup (1-\e,1]$ implies $\|1_A-\alpha_{\Sigma(d)}\|_{U^2(d)}<2\e^{1/2}$.

Conversely, if $A$ has $\VC$-dimension at most $k$, then $\|1_A-\alpha_{\Sigma(d)}\|_{U^2(d)}<\e$ implies that $\alpha_{\Sigma(d)}\in [0,\e^{1/2k})\cup (1-\e^{1/2k},1]$.
\end{fact}

By Fact \ref{fct:unihom}, point (iii) in Corollary \ref{cor:linearregvc} implies that for all $d\in \Gamma$, $\|1_A-\alpha_{\Sigma(d)}\|_{U^2(d)}<2\e^{1/2}$, so Corollary \ref{cor:linearregvc} is indeed a fully-fledged arithmetic regularity lemma.   We note that the conclusion of Corollary \ref{cor:linearregvc} is quantitatively stronger than that of Theorem \ref{thm:linearreg} because of the polynomial bound on $\ell$. It is also qualitatively stronger because, in the absence of the assumption of bounded VC-dimension, uniformity does not necessarily imply density near 0 or near 1. For example, a set $A\subseteq \F_p^n$ whose elements are chosen at random with probability $1/2$ will, with high probability, be uniform with density around $1/2$ on most atoms of a linear factor of bounded complexity. 
 
Work of Gowers on Szemer\'edi's theorem in the late 1990s  \cite{Gowers.1998, Gowers.2001}  introduced higher-order notions of uniformity, encapsulated by the $U^k$ uniformity norms that now carry his name. These were designed to handle more complex problems in additive groups than the original notion of uniformity based on the Fourier transform. This work set the stage for the development of \emph{higher-order arithmetic regularity lemmas}. In this paper, we shall be interested exclusively in the first non-linear level of this higher-order hierarchy, namely the quadratic one, which corresponds to uniformity as measured by the $U^3$ norm and relies crucially a so-called inverse theorem for this norm proved by Green and Tao in the mid-2000s \cite{Green.2008piq}. 

Green and Tao also introduced the terminology of \emph{quadratic factors}, which, roughly speaking, correspond to partitions of $\F_p^n$ into the level sets (known as \emph{quadratic atoms}) of a bounded number of linear and quadratic forms. The \emph{complexity} of a quadratic factor, which is closely related to the partition size, is measured by two parameters, the number of linear forms $\ell$ and the number of quadratic forms $q$. We refer the reader to Section \ref{subsec:linquadfac} for a brief summary of the terminology of quadratic factors, to \cite{Terry.2021d} for a more comprehensive overview, and to \cite{Green.2007} for a more in-depth introduction.  For what follows, it suffices to know a quadratic factor of complexity $(\ell,q)$ takes the form of a pair $\calB=(\calL,\calQ)$, where $\calL$ is a linear factor of complexity $\ell$, and $\calQ$ is a ``purely quadratic" factor of complexity $q$, which comes equipped with a notion of rank. Such a quadratic factor $\calB$ partitions the group $\F_p^n$ into \emph{atoms}, each of which is denoted $B(a)$ for some $a\in \F_p^{\ell+q}$.

A quadratic arithmetic regularity lemma then states most naturally, in view of Theorem \ref{thm:linearreg}, that given any set $A\subseteq \F_p^n$, there exists a bounded number of linear and quadratic forms such that on almost all atoms of the resulting quadratic factor, $A$ is uniform as measured by the $U^3$ norm. In order to express Green and Tao's result \cite{Green.2007} in this way, we introduced a local version of the $U^3$ norm in \cite{Terry.2021d}, which allowed us to rephrase the quadratic arithmetic regularity lemma as follows (see \cite[Proposition 4.2]{Terry.2021d}).
 
 \begin{theorem}[$U^3$ arithmetic regularity \cite{Green.2007}, \cite{Terry.2021d}]\label{thm:localu3}
 For all primes $p>2$ and $\e\in (0,1)$, there is a polynomial growth function $\rho_1=\rho_1(\e,p):\mathbb{R}^+\rightarrow \mathbb{R}^+$ such that for all growth functions $\rho\geq\rho_1$ (pointwise),  there is a positive integer $D=D(\e,\rho)$ such that the following holds for all sufficiently large $n$.
 
 Given $A\subseteq \F_p^n$, there exist integers $\ell,q\geq 1$ and a quadratic factor $\calB=(\calL,\calQ)$ on $\F_p^n$ of complexity $(\ell, q)$ together with a set $\Gamma\subseteq \mathbb{F}_p^{3\ell+6q}$ with the property that
\begin{enumerate}[label=\normalfont(\roman*)]
\item $\ell,q\leq D$;
\item  $\calB$ has rank at least $\rho(\ell+q)$;
\item $|\Gamma|\geq (1-\e)p^{3\ell+6q}$;
\item for all $d=(a_1,a_2,a_3,b_{12},b_{13},b_{23})\in \Gamma$, $\|1_A-\alpha_{B(\Sigma(d))}\|_{U^3(d)}<\e$, where $\alpha_{B(\Sigma(d))}$ is the density of $A$ on the atom $B(\Sigma(d))=B(a_1+a_2+a_3+2(0b_{12})+2(0b_{13})+2(0b_{23}))$.
\end{enumerate}
 \end{theorem}
 
The definitions of several terms appearing in  Theorem \ref{thm:localu3} are not strictly necessary for this paper.  For the purposes of this introduction, it suffices to remark that  the atoms of the form $B(\Sigma(d))$ appearing in (iv) are considered \emph{uniform with respect to $A$}, and the bound in (iii) ensures almost all atoms of $\calB$ are uniform in this sense.
 
 For reasonable rank functions $\rho$, the proof in \cite{Green.2007} produced a tower-type bound on $D$ in terms of the inverse of the uniformity parameter $\e$ (see e.g. the appendix to \cite{Gladkova.2025}). Recent work of Gladkova \cite{Gladkova.2025} also shows that the bound on the complexity $\ell$ of the linear part of the factor obtained in this proof of Theorem \ref{thm:localu3} cannot be improved, in analogy with a result of Moshkovitz and Shapira \cite{Moshkovitz.2019} for 3-uniform hypergraphs.

In \cite{Terry.2021d}\footnote{\label{note1}See \cite[Theorem 1.9]{Terry.2021d}, which was first proved in an early version of \cite{Terry.2021a} as \cite[Theorem 1.8]{Terry.2021av2}.}, we went on to prove a stronger version of Theorem \ref{thm:localu3} under the assumption that the subset in question has bounded $\VC_2$-dimension. This ternary analogue of $\VC$-dimension was first defined by Shelah in \cite{Shelah.2014}. It was shown to be related to hypergraph regularity and other combinatorial problems in \cite{Terry.2018, Chernikov.2020, Terry.2021b}, as well as bilinear forms over vector spaces in \cite{Chernikov.2019b, Hempel.2016}. In \cite{Terry.2021d}, the authors studied this notion in the context of subsets of abelian groups.

\begin{definition}[$\VC_2$-dimension of a subset of a group]\label{def:vc2dim}
Given an integer $k\geq 1$, an abelian group $G$ and a subset $A\subseteq G$, we say that $A$ \emph{has $\VC_2$-dimension at least $k$} (or \emph{has $k$-$\IP_2$}) if there exist 
$$
\{a_i:i \in [k]\}\cup \{b_{j}:j\in [k]\}\cup \{c_S: i\in [k]\times [k]\}\subseteq G
$$
 such that $a_i+b_j+ c_S \in A$ if and only if $(i,j)\in S$.

The \emph{$\VC_2$-dimension of $A$}, denoted $\VC_2( A)$, is defined to be the largest integer $k\geq 1$ such that $A$ has $\VC_2$-dimension at least $k$, when such a $k$ exists, and is otherwise defined to be $0$.
\end{definition}

We are now able to state the main result of \cite{Terry.2021d}, which is a higher order version of Theorem \ref{thm:vc}.\footnote{See footnote \ref{note1}.}

\begin{theorem}[Structure theorem for sets of bounded $\VC_2$-dimension \cite{Terry.2021d}]\label{thm:vc2}
Let $p>2$ be a fixed prime. For all integers $k\geq 0$ and $\mu\in (0,1)$, there exists a polynomial growth function $\omega_0=\omega_0(p,k,\mu):\mathbb{R}^+\rightarrow \mathbb{R}^+$ such that for all $\omega\geq\omega_0$ (pointwise) there exists a positive integer $C=C(p,k,\mu,\omega)$ such that the following holds for all sufficiently large $n$.     

Given $A\subseteq \F_p^n$ with $\VC_2(A)\leq k$, there exist integers  $1\leq \ell,q\leq C$, a quadratic factor $\calB=(\calL,\calQ)$ on $\F_p^n$ of complexity $(\ell, q)$ and rank at least $\omega(\ell+q)$,  and a union $Y$ of atoms of $\calB$ so that $|A\Delta Y|\leq \mu |G|$.
\end{theorem} 

Combining \cite[Proposition 3.20]{Terry.2021d} and \cite[Proposition 4.3]{Terry.2021d}, we have the following higher-order analogue of Fact \ref{fct:unihom}, which states that under the assumption of bounded $\VC_2$-dimension, having density near 0 or 1 is approximately equivalent to satisfying our higher-order notion of uniformity.

\begin{fact}[Local homogeneity and local $U^3$ uniformity are approximately equivalent \cite{Terry.2021d}]\label{fct:unihom2}
Let $p>2$ be a fixed prime. For all integers $k\geq 0$ and $\e\in (0,1)$, there exists a polynomial growth function $\rho_0=\rho_0(\e,k):\mathbb{R}^+\rightarrow \mathbb{R}^+$ such that for all growth functions $\rho\geq \rho_0$ (pointwise), all integers $\ell,q\geq 1$ and all $n > \ell+q$ the following holds.

Let $\calB=(\calL,\calQ)$ be a quadratic factor on $\F_p^n$ of complexity $(\ell,q)$ and rank at least $\rho(\ell+q)$, let $d=(a_1,a_2,a_3,b_{12},b_{13},b_{23})\in \F_p^{\ell+q}\times \F_p^{\ell+q}\times \F_p^{\ell+q}\times\F_p^q\times \F_p^q\times \F_p^q$, and $A\subseteq \F_p^n$ has density $\alpha_{\Sigma(d)}$ on the atom $B(\Sigma(d))=\Sigma(d)=a_1+a_2+a_3+ 2(0b_{12})+2(0b_{13})+2(0b_{23})$.

Then $\alpha_{\Sigma(d)}\in [0,\e)\cup (1-\e,1]$ implies $\|1_A-\alpha_{\Sigma(d)}\|_{U^3(d)}<2 \e^{1/8}$.

Conversely, if $A$ has $\VC_2$-dimension at most $k$, then $\|1_A-\alpha_{\Sigma(d)}\|_{U^3(d)}<\e$ implies that $\alpha_{\Sigma(d)}\in [0,\e^{1/(k^2 2^{k^2})})\cup (1-\e^{1/(k^2 2^{k^2})},1]$.
\end{fact}

We note that if a set $A$ does not have bounded $\VC_2$-dimension, then it is possible for an atom $B(\Sigma(d))$ to be uniform with respect to $A$ in the sense of Theorem \ref{thm:localu3} while having density bounded away from $0$ and $1$. Indeed, a set whose elements are chosen randomly with probability $1/2$ will be uniform with density close to $1/2$ on most atoms of a bounded complexity factor,  with high probability.  

With Fact \ref{fct:unihom2} in hand, we see that Theorem \ref{thm:vc2} gives rise to a strengthened version of Theorem \ref{thm:localu3} in which the overwhelming majority of atoms are not just uniform, but have density near 0 or 1. 

\begin{corollary}[$U^3$ regularity for sets of bounded $\VC_2$-dimension \cite{Terry.2021d}]\label{cor:vc2}
Let $p>2$ be a fixed prime. For all integers $k\geq 0$ and $\mu\in (0,1)$, there exists a polynomial growth function $\omega_0=\omega_0(p,k,\mu):\mathbb{R}^+\rightarrow \mathbb{R}^+$ such that for all $\omega\geq\omega_0$ (pointwise) there exists $C=C(p,k,\mu,\omega)$  such that the following holds for all sufficiently large $n$. 

Given $A\subseteq \F_p^n$ with $\VC_2(A)\leq k$, there exist integers $\ell,q\geq 1$ and a quadratic factor $\calB=(\calL,\calQ)$ on $\F_p^n$ of complexity $(\ell, q)$ together with a set $\Gamma\subseteq \mathbb{F}_p^{3\ell+6q}$ with the property that
\begin{enumerate}[label=\normalfont(\roman*)]
\item $\ell, q\leq C$;
\item $\calB$ has rank at least $\omega(\ell+q)$;
\item $|\Gamma|\geq (1-\e)p^{3\ell+6q}$; 
\item for all $d=(a_1,a_2,a_3,b_{12},b_{13},b_{23})\in \Gamma$, $|A\cap B(\Sigma(d))|/|B(\Sigma(d))|\in [0,\e)\cup (1-\e,1]$. 
\end{enumerate}
\end{corollary}

Note again that by Fact \ref{fct:unihom2}, for all $d=(a_1,a_2,a_3,b_{12},b_{13},b_{23})\in \Gamma$, it follows from (iv) in Corollary \ref{cor:vc2} that we also have $\|1_A-\alpha_{B(\Sigma(d))}\|_{U^3(d)}<2\e^{1/8}$, so Corollary \ref{cor:vc2} is a true higher-order regularity statement.

In light of Corollary \ref{cor:linearregvc}, one would expect the bounds on the complexity of the quadratic factor obtained under the additional assumption of bounded $\VC_2$-dimension to be rather stronger than in the case of arbitrary subsets of $\F_p^n$. However, the initial approach taken in \cite{Terry.2021d} to prove Theorem \ref{thm:vc2} used Theorem \ref{thm:localu3} as a starting point, and thus gave rise to tower-type bounds (in $\mu^{-1}$) on both $\ell$ and $q$.

The main purpose of the present paper is to improve the bound on the quadratic complexity $q$ to a logarithm in $\mu^{-1}$, meaning that the partition of the group generated by the purely quadratic portion of the factor has size at most polynomial in $\mu^{-1}$. More precisely, we shall prove the following higher arity analogue of Theorem \ref{thm:vc}.

\begin{theorem}[Structure theorem for sets of bounded $\VC_2$-dimension with improved quadratic complexity]\label{thm:vc2main}
Let $p>2$ be a fixed prime. For all integers $k\geq 0$ and $\mu\in (0,1)$, there exists a polynomial growth function $\omega_0=\omega_0(p,k,\mu):\mathbb{R}^+\rightarrow \mathbb{R}^+$ such that for all $\omega\geq\omega_0$ (pointwise) there exists $C=C(p,k,\mu, \omega)$ such that the following holds for all sufficiently large $n$.  

Given $A\subseteq \F_p^n$ with $\VC_2(A)\leq k$, there exist integers $\ell,q\geq 1$ and a quadratic factor $\calB=(\calL,\calQ)$ on $\F_p^n$ of complexity $(\ell,q)$, such that
\begin{enumerate}[label=\normalfont(\roman*)]
\item $\ell \leq C$ and  $q\leq \log_p(\mu^{-k-o_{k,q}(1)})$, where $o_{k,p}(1)$ tends to $0$ with $\mu$, at a rate depending on $k$ and $p$;
\item  $\calB$ has rank at least $\omega(\ell+q)$; 
\item for some union $Y$ of atoms of $\calB$, $|A\Delta Y|\leq \mu |G|$.
\end{enumerate}
\end{theorem}

From this we will deduce a similar quantitative improvement of Theorem \ref{thm:localu3} in the setting of bounded $\VC_2$-dimension.

\begin{corollary}[$U^3$ regularity for sets of bounded $\VC_2$-dimension with improved quadratic complexity]\label{cor:betterbound}
Let $p>2$ be a fixed prime.  For all integers $k\geq 0$ and $\e\in (0,1)$,  there is a polynomial growth function $\rho_0=\rho_0(p,k,\e):\mathbb{R}^+\rightarrow \mathbb{R}^+$ such that for all growth functions $\rho\geq\rho_0$ (pointwise), there exist $D=D(p,k,\e,\rho)$ such that the following holds for all sufficiently large $n$.  
 
Given $A\subseteq \F_p^n$ with $\VC_2(A)\leq k$, there exist integers $\ell,q\geq 1$ and a quadratic factor $\calB=(\calL,\calQ)$ on $\F_p^n$ of complexity $(\ell, q)$, together with a set $\Gamma\subseteq \mathbb{F}_p^{3\ell+6q}$ with the property that
\begin{enumerate}[label=\normalfont(\roman*)]
\item $\ell \leq D$ and $q\leq \log_p(\e^{-16k-o_{k,p}(1)})$, where $o_{k,p}(1)$ tends to $0$ as $\e$ tends to $0$ at a rate depending on $p$ and $k$;
\item  $\calB$ has rank at least $\rho(\ell+q)$;
\item $|\Gamma|\geq (1-\e)p^{3\ell+6q}$;
\item for all $d=(a_1,a_2,a_3,b_{12},b_{13},b_{23})\in \Gamma$, $\|1_A-\alpha_{B(\Sigma(d))}\|_{U^3(d)}<\e$, and moreover, $|A\cap B(\Sigma(d))|/|B(\Sigma(d))|\in [0,\e^8/2)\cup (1-\e^8/2,1]$.
\end{enumerate}
\end{corollary}

In fact, our proof shows that any quadratic factor $\calB$ that is sufficiently uniform with respect to $A$ (in the sense of Theorem \ref{thm:localu3}) can be turned into one whose quadratic complexity is efficiently bounded (see Remark \ref{rem:strength}). 

Finally, we give an example of a subset of $\F_p^n$ with the property that any quadratic arithmetic regularity decomposition has quadratic complexity that grows at least like a power of $\e^{-1}$, showing that Corollary \ref{cor:betterbound} does not hold in general.
 
\begin{theorem}\label{thm:lowerbound}
Let $p>2$ be a fixed prime. There exist  $\e\in (0,1)$ and a polynomial growth function $\rho:\mathbb{R}^+\rightarrow \mathbb{R}^+$ such that for all sufficiently large odd integers $n$,  there exists $A\subseteq\F_p^n$ so that the following holds.

Let $\ell,q\geq 1$ be integers, and let $\calB=(\calL,\calQ)$ be a quadratic factor on $\F_p^n$ of complexity $(\ell,q)$ and rank at least $\rho(\ell+q)$.  Suppose that $\Gamma\subseteq \F_p^{3\ell+6q}$ satisfies
\begin{enumerate}[label=\normalfont(\roman*)] 
\item $|\Gamma|\geq (1-\e)p^{3\ell+6q}$; and
\item for all $d\in \Gamma$, $\|1_A-\alpha_{B(\Sigma(d))}\|_{U^3(d)}<\e$.
\end{enumerate}
Then $q \geq (1-\e)\e^{-1/64}$.
\end{theorem}
 
It seems plausible that there is a construction that forces the quadratic complexity to be even larger (possibly of tower-type), but this is outside the scope of the current paper.

Corollary \ref{cor:betterbound} is analogous to work of the first author in the context of hypergraphs \cite{Terry.2022}. In particular, \cite{Terry.2022} shows that $3$-uniform hypergraphs of bounded $\VC_2$-dimension satisfy a hypergraph regularity lemma in which the complexity of the partition of pairs of vertices is bounded by a polynomial in the regularity parameter. It is possible to translate Corollary \ref{cor:betterbound} into a regularity lemma for the ternary sum (hyper)graph of the set $A$ (see \cite[Appendix A]{Terry.2021a}). The complexity of the partition of pairs of vertices of the resulting decomposition will be $p^q\leq \mu^{-O_k(1)}$, matching the bound in \cite{Terry.2022}.  

The proof of Theorem \ref{thm:vc2main} is reminiscent of elements of \cite{Terry.2021a}, in that we use Theorem \ref{thm:vc2} to generate a certain ``reduced" set associated to a set $A\subseteq \F_p^n$ of bounded $\VC_2$-dimension in $\mathbb{F}_p^n$. This reduced set comprises the set of labels of those quadratic atoms arising from Theorem \ref{thm:vc2} on which $A$ has density close to 1, and thus lives in $\F_p^{\ell+q}$ with $\ell+q$ much smaller than $n$ (also known as ``configuration space"). The reduced set can be shown to have approximately bounded $\VC$-dimension relative to a special subgroup, namely the one generating the linear part of the factor in Theorem \ref{thm:vc2}. To this reduced set we apply a strengthened version of Theorem \ref{thm:vc} in the smaller group $\F_p^{\ell+q}$, the result of which can be pulled back to the original group $\F_p^n$ to obtain Theorem \ref{cor:betterbound}.  The difference between this argument and that in \cite{Terry.2021a} is that here we apply a structure theorem related to $\VC$-dimension in configuration space, whereas in \cite{Terry.2021a} we applied a structure theorem related to stability. 

\textbf{Outline of the paper.} We begin in Section \ref{sec:lochaussler} by recalling the strategy used by Alon, Fox and Zhao to prove Theorem \ref{thm:vc} before adapting it to work relative to a fixed subspace. In Section \ref{sec:approxAFZ}, we further generalise this result to work in the setting of 3-colorings where one of the colour classes is small, approximating the situation where we consider a set and its complement as a 2-coloring of the space. In Section \ref{sec:suffcondip2}, after giving some basic background on quadratic factors, we state a sufficient condition (from \cite{Terry.2021a}) for the existence of $\IP_2$ based on the presence of $\IP$ with one side in the subspace generating the linear part of factor given by Theorem \ref{thm:vc2}. Finally, in Section \ref{sec:betterproof} we obtain the improvement on the quadratic complexity of the factor promised in Theorem \ref{thm:vc2main} by combining the earlier ingredients with some pull-back arguments from \cite{Terry.2021a}. Section \ref{sec:lowerbound} contains the lower bound construction that proves Theorem \ref{thm:lowerbound}.

\textbf{Notation.} 
For real numbers $a$, $b$ and $c$, we will occasionally use the notation $a= b \pm c$ to mean $|a-b|\leq c$. Throughout, given $v\in \F_p^n$ and $1\leq i\leq n$, we let $v_i$ denote the $i$th coordinate of the vector $v$.


\textbf{Acknowledgements.}
The first author is partially supported by NSF CAREER Award DMS-2115518 and a Sloan Research Fellowship. The second author is supported by an Open Fellowship from the UK Engineering and Physical Sciences Research Council (EP/Z53352X/1).

\section{A localized version of the Alon-Fox-Zhao theorem}\label{sec:lochaussler}

In this section we prove a structure theorem for sets which have bounded $\VC$-dimension relative to a subgroup. The proof of this result will later serve as a template for a more complicated version needed in the proof of our main theorem.  We begin by sketching the argument used by Alon, Fox and Zhao \cite{Alon.2018is} to prove Theorem \ref{thm:vc}, as we shall need several of the ingredients in due course. A crucial role in this argument (and indeed, in model theoretic work on this problem) is played by the so-called \emph{stabiliser} of a set. 

\begin{definition}[Stabilizer of a set]\label{def:stabiliser}
Given $A\subseteq G= \F_p^n$ and $\rho>0$, define 
\[\Stab_\rho(A):=\{x\in G: |A\Delta (A+x)|\leq \rho|G|\}.\]
\end{definition}

Alon, Fox and Zhao showed that in the setting of finite abelian groups, a subset $A$ of bounded $\VC$-dimension must have a large stabilizer. The idea of finding large stabilizers in the context of bounded $\VC$-dimension had previously appeared in model theory. For example, this idea is implicitly present in \cite[Corollary 4.3]{Hrushovski.2008}.  

In the setting of finite abelian groups, this is achieved using a foundational result in $\VC$-theory known as Haussler's packing lemma \cite{Haussler.1995} (see also \cite{Matousek.2002}, for example). Before we state it, recall that two subsets $F,F'\subseteq V$ of a set $V$ are said to be \emph{$\delta$-separated} for $\delta>0$ if $|F\cap F'|\geq \delta|V|$.  Given a collection of subsets of $V$, $\calF\subseteq \calP(V)$, we say $\calF$ is \emph{$\delta$-separated}, if all of the sets in $\calF$ are pairwise $\delta$-separated.

\begin{theorem}[Haussler's packing lemma]\label{thm:haussler}
Let $\calF\subseteq \calP(V)$ be a collection of subsets of a set $V$ and assume that $k\geq 0$. Suppose that $\calF$ has $\VC$-dimension at most $k$.  Then any $\delta$-separated subcollection of $\calF$ has size at most $(30/\delta)^k$.
\end{theorem}

The first step in the Alon-Fox-Zhao strategy is to deduce deduce from this that if $A$ has bounded $\VC$-dimension, then there is a bounded number of  ``typical" translates of $A$ with the property that every translate of $A$ looks approximately like one of the typical ones.

\begin{proposition}[Representation lemma]\label{prop:gphaussler}
Suppose that $A\subseteq G=\F_p^n$ has $\VC$-dimension at most $k$. Then there exists an integer $m\leq (30/\delta)^k$ and elements $g_1,\ldots, g_m\in \F_p^n$ such that for all $x\in G$, there exists $i\in [m]$ such that 
\[|(A-x)\Delta (A-g_i)|\leq \e|G|.\]
\end{proposition}

As was shown in \cite{Alon.2018is}, one can deduce from Proposition \ref{prop:gphaussler} that the stabilizer of $A$ is of positive density in $G$.  Indeed, the conclusion of Proposition \ref{prop:gphaussler} implies that for each $x\in G$, there exists some $i\in [m]$ such that $x-g_i\in \Stab_\e(A)$. Consequently, $G\subseteq  \bigcup_{i=1}^m(\Stab_{\e}(A)+g_i)$, yielding that $|\Stab_{\e}(A)|\geq |G|/m$. 

The second step in the Alon-Fox-Zhao strategy is to identify  a subgroup inside the dense stabiliser obtained in Step 1. This is done using standard additive-combinatorial tools, including Pl\"unnecke's inequality and the Freiman-Ruzsa theorem (see \cite[Lemma 2.4]{Alon.2018is}), where the parameters are chosen with the particular application in mind. 

\begin{proposition}[Sumsets of dense sets contain a subgroup]\label{prop:stabafz}
Let $p$ be a prime, let $\delta\in (0,1/2)$, and let $C>0$ a constant. Suppose that $B\subseteq G=\F_p^n$ satisfies $|B|\geq \delta^C|G|$.  Then there is a positive integer $\ell \leq \delta^{-o_{C,p}(1)}$ and a subgroup $H\leqslant G$ satisfying $H\subseteq 2\ell B-2\ell B$ and $|H|\geq \delta^{o_{C,p}(1)}|B|$, where $o_{C,p}(1)\rightarrow 0$ as $\delta\rightarrow 0$ at a rate depending on $C$ and $p$.

Specifically, the $o_{k,p}(1)$ term can be taken to have the form $\phi(p,C,\delta)=\gamma_{p,C}(\log(\delta^{-1}))^{-1/5}$, where $\gamma_{p,C}$ is a constant depending on $C$ and $p$.
\end{proposition}

We remark that the explicit description of the $o_{C,p}(1)$ term in Proposition \ref{prop:stabafz} will not be used until Section \ref{sec:approxAFZ}.  

The third and final step in  the Alon-Fox-Zhao strategy is to observe that if there is a subgroup $H$ contained in a stabilizer of a set $A$, then $A$ is  is approximately equal to a union of cosets of $H$. This is used in the proof of \cite[Lemma 2.4]{Alon.2018is}, and also appears explicitly in \cite[Lemma 8.2]{Conant.2020}. 

\begin{lemma}\label{lem:stabsubgp}
Suppose that $A\subseteq G=\F_p^n$ and $H \leqslant G$ is a subgroup such that $H\subseteq \Stab_{\e}(A)$.  Then there is a union $Y$ of cosets of $H$ so that $|A\Delta Y|\leq \e|G|$. 
\end{lemma}

Combining Proposition \ref{prop:gphaussler}, Lemma \ref{lem:stabsubgp} and Proposition \ref{prop:stabafz} yields Theorem \ref{thm:vc}.

Our goal in this section is to prove a version of Proposition \ref{prop:gphaussler} relative to a fixed subspace, before deducing a structural result analogous to Theorem \ref{thm:vc} for sets with bounded $\VC$-dimension relative to a subspace, which we now define. 

\begin{definition}[$\VC$-dimension relative to a subgroup]\label{def:vcsub}
Let $k\geq 1$ be an integer, let $G$ be a finite abelian group, let $H\leqslant G$ be a subgroup, and let $A\subseteq G$.  We say that $A$ has \emph{$\VC$-dimension at least $k$ relative to $H$} if there exist $x_1,\ldots, x_k\in G$ and for each $S\subseteq [k]$, some $y_S\in H$ so that $x_i+y_S\in A$ if and only if $i\in S$.

The \emph{$\VC$-dimension of $A$ relative to $H$} is defined to be the largest integer $k\geq 1$ for which the $\VC$-dimension of $A$ relative to $H$ is at least $k$, if such a $k$ exists, and is otherwise defined to be $0$.
\end{definition}

Observe that in the notation of Definition \ref{def:vcsub}, the $\VC$-dimension of $A$ relative to $H$ is simply the VC-dimension of the set system
$$
\calF=\{A-h: h\in H\}
$$
on the ground set $G$. We note that Definition \ref{def:vcsub} is the same as the notion $\VC_H(A)$ used in \cite{Conant.2025}, and the dual of $\dim_{\VC}(A|H)$ used in \cite{Sisask.2018}. For further discussion on this and related notions of $\VC$-dimension, see \cite[Section 4]{Sisask.2018} and \cite[Appendix A]{Conant.2025}.

We may now deduce the following analogue of Proposition \ref{prop:gphaussler}.

\begin{proposition}[Representation lemma relative to a subspace]\label{prop:hausslerlocal}
Let $p$ be a prime, let $k\geq 0$ be an integer, and let $\delta\in (0,1)$. Let $H\leqslant G=\F_p^n$ be a subgroup, and let $A\subseteq G$  be a subset whose $\VC$-dimension relative to $H$ is at most $k$.
Then there exists an integer $m\leq (30/\delta)^{k}$ and elements $h_1,\ldots, h_m\in H$ such that for all $x\in H$, there is $i\in [m]$ satisfying
\[|(A-x)\Delta (A-h_i)|\leq \delta|G|.\]
\end{proposition}

\begin{proof}
By definition,  the set system $\calF=\{A-x: x\in H\}$ on $G$ has $\VC$-dimension at most $k$.  Let $X=\{h_1,\ldots, h_m\}$ be a maximal $\delta$-separated subcollection of $\calF$.  By Haussler's Lemma (Theorem \ref{thm:haussler}), $m\leq (30/\delta)^{k}$. By construction, we know that for all $x\in H$, there is $i\in [m]$ such that 
\[|(A-x)\Delta (A-h_i)|\leq \delta |G|,\]
which finishes the proof.
 \end{proof}
 
 We next prove that given a set $A$ whose VC-dimension relative to a fixed subgroup $H$ is bounded, we can find a translate of $A$ whose stabilizer has large intersection with the subgroup $H$.  We note this also follows from  \cite[Proposition 2.10]{Conant.2025}, which uses a similar strategy.  We include another proof here, as we will use it as a model for a more complicated argument in the next section.

\begin{lemma}\label{lem:mainvcrel}
Let $p$ be a prime, let $k\geq 0$ be an integer, and let $\delta\in (0,1)$. Let $H\leqslant G=\F_p^n$, and let $A\subseteq G$  be a subset whose $\VC$-dimension relative to $H$ is at most $k$. Then 
$$
|\Stab_{\delta}(A)\cap H|\geq (\delta/30)^{k}|H|. 
$$
\end{lemma}
\begin{proof} 
By Proposition \ref{prop:hausslerlocal}, there is an integer $m\leq (30/\delta)^{k}$ and elements $h_1,\ldots, h_m\in H$ so that for all $x\in H$, there is $i\in [m]$ so that $|(A-x)\Delta (A-h_i)|\leq \delta |G|$.  By the pigeonhole principle, there exists some $i\in [m]$ such that the set
$$
S=\{x\in H: |(A-x)\Delta (A-h_i)|\leq \delta|G|\}
$$
has size at least $|H|/m$.

We claim that for all $h\in S-h_i$, $|(A-h)\Delta A|\leq  \delta |G|$.  To this end, fix $h\in S-h_i$. Then we can write $h=s -h_i$ for  some $s\in S$, and consequently,
\begin{align*}
|(A-h)\Delta A|=|(A-s+h_i)\Delta A|=|(A-s)\Delta (A-h_i)|&\leq \delta |G|,
\end{align*}
where the last inequality holds because $s\in S$.  Hence $h\in \Stab_{\delta}(A)\cap H$. We have now shown $S-h_i\subseteq \Stab_\delta(A)\cap H$, hence $|\Stab_{\delta}(A)\cap H|\geq |S|\geq |H|/m\geq (\delta/30)^{k}|H|$, as desired.
\end{proof}

Note that the preceding two lemmas do not require the group to be abelian, a fact exploited in \cite{Conant.2025}.
Combining Lemma \ref{lem:mainvcrel} with Proposition \ref{prop:stabafz}, we obtain a structure theorem for sets with bounded VC-dimension relative to a subgroup. We shall not make use of this result, but the proof serves as a warm-up for the proof of Theorem \ref{thm:errorstructure} later on.

\begin{proposition}[Structure theorem for sets of bounded $\VC$-dimension relative to a subspace]\label{prop:strucrelvc}
Let $p$ be a prime, let $k\geq 0$ be an integer, and let $0<\delta<1/60$. Let $H\leqslant G=\F_p^n$ be a subgroup, and let $A\subseteq G$  be a subset whose $\VC$-dimension relative to $H$ is at most $k$.  Then there exists a subgroup $H'\leqslant H$ such that  
\begin{enumerate}[label=\normalfont(\roman*)]
\item $H'$ has index at most $\delta^{-k-o_{k,p}(1)}$ in $H$,
\item for all $h\in H'$, $|A\Delta (A+h)|\leq 4\delta^{1-o_{k,p}(1)} |G|$, and
\item there is a union $Y$ of cosets of $H'$ such that $|A\Delta Y|\leq 4\delta^{1-o_{k,p}(1)} |G|$,
\end{enumerate}
where $o_{k,p}(1)\rightarrow 0$ as $\delta\rightarrow 0$ at a rate which depends on $k$ and $p$.
\end{proposition}
\begin{proof}
By Lemma \ref{lem:mainvcrel},
$$
|\Stab_{\delta}(A)\cap H|\geq (\delta/30)^{k}|H|. 
$$
Define $B=\Stab_{\delta}(A)\cap H $ and let $D>0$ be a constant such that $|B|=\delta^D|H|$.
Since $|B|\geq (\delta/30)^{k}|H|$, we have 
$$
D\leq k\Big(1+\frac{\log 30}{\log (\delta^{-1})}\Big)=k+o_k(1),
$$
where $o_k(1)=k\log(30)/\log(\delta)$. Note that since $\delta<1/60$, this implies $D\leq 2k$.   Apply Proposition \ref{prop:stabafz} to $B$ (considered as a subset of the group $H$), and parameters $\delta$ and $C=2k$ to obtain $\ell \leq \delta^{-o_{C,p}(1)}$ and $H'\leqslant H$ satisfying $H'\subseteq 2\ell B-2\ell B$ and $|H'|\geq \delta^{o_{C,p}(1)}|B|$, where $o_{C,p}(1)$ is the function from Proposition \ref{prop:stabafz}  tending to $0$ with $\delta$, at a rate depending on $C$ and $p$ (and thus at a rate depending on $k$ and $p$ since $C=2k$).   The lower bound on $|H'|$ implies that the index of $H'$ in $H$ is at most 
$$
\frac{|H|}{\delta^{o_{C,p}(1)}|B|}= \frac{1}{\delta^{D+o_{k,p}(1)}}\geq \delta^{-k-o_k(1)-o_{C,p}(1)},
$$
where the last equality use that $D\leq k+o_k(1)$.  Since $H'\subseteq 2\ell B-2\ell B$ and $B\subseteq \Stab_{\delta}(A)$, we have by the triangle inequality that for all $x\in H'$,
\begin{align}\label{adelta}
|A\Delta (A +x)|\leq  4\ell \delta |G|\leq 4\delta^{1-o_{C,p}(1)} |G|\leq 4\delta^{1-o_k(1)-o_{C,p}(1)} |G|.
\end{align}
We have now shown (i) and (ii) hold with the error function $o_{k,p}(1)=o_{C,p}(1)+o_{k}(1)$.  Combining (ii) and Lemma \ref{lem:stabsubgp} implies that there is a union $Y$ of cosets of $H'$ so that $|(A-g)\Delta Y|\leq 4\delta^{1-o_{k,p}(1)}  |G|$, yielding (iii).
\end{proof}

\section{An approximate version of the local Alon-Fox-Zhao theorem}\label{sec:approxAFZ}

 Theorem \ref{thm:vc2} will, given a subset $A\subseteq \F_p^n$ of small $\VC_2$-dimension, output a smaller group of the form $\F_p^m$ (for some $m$ that depends only on the parameters $p$, $k$ and $\e$, but which is crucially independent of the starting dimension $n$), along with a distinguished subspace of $\F_p^m$ and three types of elements in $\F_p^m$, corresponding to ``dense", ``sparse", and ``intermediate" atoms of the quadratic decomposition.  We refer the reader to the next section for details on this construction. It suffices to say here that for this reason, we need to work with groups equipped with a distinguished partition into three parts.  

\begin{definition}[3-coloring of a group]
A \emph{$3$-coloring of a group $G$} is a tuple $(G;A_0,A_1,A_2)$ where  $A_0\cup A_1\cup A_2$ is a partition of $G$.  
\end{definition}
We next define analogues of VC-dimension relative to a subgroup for the setting of $3$-colorings.  This definition first appeared in \cite{Terry.2021a} (see Definition 4.20).

\begin{definition}[$A_1|A_0$-copy of $k$-$\IP$ with right side in $H$]
Let $H\leqslant G=\F_p^n$, and let $P=(G;A_0,A_1, A_2)$ be a $3$-coloring.  We say that  \emph{$A_1|A_0$ has $\VC$-dimension at least $k$ relative to $H$} if there are $v_1,\ldots, v_k\in G$ and for each $S\subseteq [k]$ an element $w_S\in H$ such that $i\in S$ implies $v_i+w_S\in A_1$ and $i\notin S$ implies $v_i+w_S\in A_0$.  

We then say \emph{$A_1|A_0$ has $\VC$-dimension at most $k$ relative to $H$} if $A_1|A_0$ does not have $\VC$-dimension at least $k+1$ relative to $H$.
\end{definition}

We now prove an analogue of Proposition \ref{prop:hausslerlocal} for a special $3$-colored setting in which the class $A_2$ is very small.

\begin{proposition}[Representation lemma for 3-colorings relative to a subspace]\label{prop:hausslerthree}
Let $p$ be a prime, let $k\geq 0$ be an integer, let $0<\delta<1/4$, and let $0< \e\leq (\delta/120)^{k+1}$.  Let $H \leqslant G=\F_p^n$, and let $P=(G;A_0,A_1,A_2)$ be a $3$-coloring with $|A_2|\leq \e|G|$.  

Suppose that $A_1|A_0$-has $\VC$-dimension at most $k$ relative to $H$. Then there is an integer $m \leq \lceil(120/\delta)^{k}\rceil$ and elements $h_1,\ldots, h_m\in H$ so that for all $h\in H$, there is $i\in [m]$ so that for each $u\in \{0,1\}$, 
$$
|(A_u-h)\Delta (A_u-h_i)|\leq \delta |G|.
$$
\end{proposition}
\begin{proof}
Let $h_1,\ldots, h_m\in H$ be a maximal collection in $H$ for which $\{A_1-h_i: i\in [m]\}$ is a $\delta/2$-separated family in $G$. We first show that $m\leq \lceil(120/\delta)^{k}\rceil$. 

Suppose towards a contradiction that $m>\lceil(120/\delta)^{k}\rceil$, so $m\geq  m_0=\lceil (\delta/120)^{-k}\rceil+1$. Define 
$$
X=G\setminus \bigcup_{i=1}^{m_0} (A_2-h_i)\text{ and }\calF:=\{(A_1-h_i)\cap X: i\in [m_0]\}.
$$
Observe that by definition, for each $i\in [m_0]$, $X\setminus (A_1-h_i)\subseteq A_0-h_i$.  Consequently, because  $A_1|A_0$ has VC-dimension at most $k$ relative to $H$, $\calF$ must have VC-dimension at most $k$, considered as a set system on $X$.

We claim that $\calF$ is $\delta/4$-separated with respect to the ground set $X$. Indeed, given $1\leq i\neq j\leq m_0$,  we have that 
\begin{align*}
|((A_1-h_i)\Delta (A_1-h_j))\cap X|&\geq |(A_1-h_i)\Delta (A_1-h_j)|-|G\setminus X|\\
&\geq |(A_1-h_i)\Delta (A_1-h_j)|-m_0\e|G|,
\end{align*}
where the last inequality holds because $|A_2|\leq \e|G|$.  Since $\{A_1-h_i: i\in [m_0]\}$ is $(\delta/2)$-separated in $G$, this implies
$$
|((A_1-h_i)\Delta (A_1-h_j))\cap X|\geq \frac{\delta|G|}{2}-m_0\e|G|\geq \frac{\delta|G|}{4}\geq \frac{\delta|X|}{4},
$$
where the first inequality holds because $m_0=\lceil (\delta/120)^{-k}\rceil+1$ and $\e$ is sufficiently small compared to $\delta$.   But now Theorem \ref{thm:haussler} implies that $\calF$ must have VC-dimension larger than $k$, a contradiction.  Thus, $m \leq  \lceil(\delta/120)^{-k}\rceil$.    

By our choice of $h_1,\dots, h_m$, we know that for all $h\in H$, there is some $i\in [m]$ such that $|(A_1-h)\Delta(A_1-h_i)|\leq \delta |G|/2$, and thus also
$$
|(A_0-h)\Delta (A_0-h_i)|\leq |(A_1-h)\Delta (A_1-h_i)|+|(A_2-h)\Delta (A_2-h_i)|\leq \frac{\delta |G|}{2}+2\e|G|\leq \delta|G|,
$$
where the inequalities use that $|A_2|\leq \e|G|$ and $\e<\delta/4$.  This finishes the proof.
\end{proof}

Following the blueprint laid out in the preceding section, we will use Proposition \ref{prop:hausslerthree}, Lemma \ref{lem:stabsubgp}, and Proposition \ref{prop:stabafz} to deduce an approximate local version of Theorem \ref{thm:vc}. We will keep closer tabs on the $o_{k,p}(1)$ terms than in the preceding section, as we will need this information later on.

\begin{theorem}\label{thm:errorstructure}
For all primes $p$ and integers $k\geq 0$, there exists $\delta_0\in (0,1)$ so that for all $0<\delta<\delta_0$ and all $0<\e\leq (\delta/120)^{k+1}$, the following holds.  Let $H\leqslant G=\F_p^n$, and let $P=(G;A_0,A_1,A_2)$ be a $3$-coloring with $|A_2|\leq \e|G|$.

Suppose that $A_1|A_0$ has $\VC$-dimension less than $k$ relative to $H$. Then exist a subspace $H'\leqslant H$  such that
\begin{enumerate}[label=\normalfont(\roman*)]
\item $H'$ has index at most $\delta^{-k-o_{k,p}(1)}$ in $H$,
\item for each $u\in \{0,1\}$ and all $h\in H'$, $|A_u\Delta(A_u-h)|\leq 4\delta^{1-o_{k,p}(1)}|G|$, and
\item for each $u\in \{0,1\}$, there is a union $Y_u$ of cosets of $H'$ in $G$, such that 
$$
|A_u\Delta Y_u|\leq 4\delta^{1-o_{k,p}(1)} |G|,
$$ 
\end{enumerate}
 where $o_{k,p}(1)\rightarrow 0$ as $\delta \rightarrow 0$ at a rate which depends on $k$ and $p$.

In particular, one can take $o_{k,p}(1)$ to have the form $\psi(p,k,\delta)=\sigma_{p,k}\log(\delta^{-1})^{-1/5}$,  where $\sigma_{p,k}$ is a constant depending on $p$ and $k$. 
\end{theorem}
\begin{proof}
Choose $\delta_0>0$ to be sufficiently small compared to $k$ and $p$. Fix $0<\delta<\delta_0$ and $0<\e<(\delta/120)^{k+1}$.   By Proposition \ref{prop:hausslerthree}, there is an  integer $m\leq \lceil (120/\delta)^{k}\rceil$ and elements $h_1,\ldots, h_m\in H$ so that for all $h\in H$,  there is some $i\in [m]$ so that for each $u\in \{0,1\}$, 
$$
|((A_u-h)\Delta (A_u-h_i))\cap H|\leq \delta |H|.
$$
By the pigeonhole principle, there is some index $i\in [m]$ so that if 
$$
S=\{x\in H: |(A_u-x)\Delta (A_u-h_i)|\leq \delta |G| \mbox{ for } u\in \{0,1\}\},
$$
then  $|S|\geq |H|/m$.  Let $B=S-h_i$.   We claim that 
\begin{align}\label{bs}
B\subseteq \Stab_{\delta}(A_0)\cap \Stab_{\delta}(A_1)\cap H.
\end{align}
To this end, fix $h\in B$.  By definition of $B$,  we can write $h=s-h_i$ for some $s\in S$. For each $u\in \{0,1\}$, we have
\begin{align*}
|(A_u-h)\Delta A_u|&=|(A_u-s+h_i)\Delta A_u|= |(A_u-s)\Delta (A_u-h_i)|\leq \delta |G|,
\end{align*}
where the last inequality holds because $s\in S$. This completes our verification of (\ref{bs}). Now let $D$ be such that $|B|=\delta^D|H|$.  Then 
$$
|B|=\delta^D|H|=|S|\geq \frac{|H|}{m}\geq   \lceil(\delta/120)^{k}\rceil|H|, 
$$
and thus $D\leq k+\psi_0(k,\delta)$, where $\psi_0(k,\delta)=\frac{k\log (120)}{\log(\delta^{-1})}$.  Since $\delta$ is sufficiently small, $D\leq 2k$.  Let $C=2k$, and let $\phi(p,C,\delta)$ be as in Proposition \ref{prop:stabafz}. Then Proposition \ref{prop:stabafz}, applied to the constants $C,\delta$ and the set $B$ implies there exist an integer $\ell\leq \delta^{-\phi(p,C,\delta)}$ and a subgroup $H'\leqslant G$ such that $H'\subseteq 2\ell B -2\ell B$ and such that
$$
|H'|\geq \delta^{\phi(p,C,\delta)}|B|\geq \delta^{D+\phi(p,C,\delta)}|H|\geq \delta^{k+\psi_0(k,\delta)+\phi(p,2k,\delta)}|H|,
$$
where the last equality uses $D\leq k+\psi_0(k,\delta)$ and the fact $C=2k$.   Since $H\subseteq 2\ell B-2\ell B$ and $B\subseteq H$,  we have that $H'\leqslant H$. Moreover, the above lower bound on $|H'|$  implies that the index of $H'$ in $H$ is at most $\delta^{-k-\psi_0(k,\delta)-\phi(p,2k,\delta)}$.   Since $H\subseteq 2\ell B-2\ell B$, the inclusion (\ref{bs}), the inequality $\ell\leq \delta^{-\phi(p,2k,\delta)}$, and the triangle inequality imply that for each $u\in \{0,1\}$ and $h\in H'$, 
$$
|A_u\Delta (A_u-h)|\leq 4\ell \delta |G|\leq 4\delta^{1-\phi(p,2k,\delta)}|G|.
$$
By Lemma \ref{lem:stabsubgp}, for each $u\in \{0,1\}$, there is a union $Y_u$ of costs of $H'$ so that 
$$
|A_u\Delta Y_u|\leq 4\delta^{1-\phi(p,2k,\delta)} |G|.
$$
Since $\psi_0(k,\delta)=k\log(120)\log(\delta^{-1})^{-1}$ and $\phi(p,2k,\delta)=\gamma_{p,2k} \log(\delta^{-1})^{-1/5}$, we have shown (i)-(iii) hold with respect to an error function of the form  $\psi(p,k,\delta)=\sigma_{p,k}\log(\delta^{-1})^{-1/5}$,  where $\sigma_{p,k}$ is a constant depending on $p$ and $k$ (e.g. let $\sigma_{k,p}=\gamma_{p,2k}+k\log(120)$). 
\end{proof}

\section{A sufficient condition for $\IP_2$ }\label{sec:suffcondip2}

In this section we set up and study the reduced colorings in configuration space that we shall be interested in. We will first recall the basic notions around linear and quadratic factors that we will need, and then develop the necessary preliminaries to state a result proved in \cite{Terry.2021a}, which says that a certain reduced set associated to a set of bounded $\VC_2$-dimension has bounded $\VC$-dimension relative to a subspace arising from the linear portion of the factor (Theorem \ref{thm:extract3}).

\subsection{Linear and quadratic factors}\label{subsec:linquadfac}

We begin by defining linear factors. 

\begin{definition}[Linear factor]\label{def:linfac}
Given an integer $\ell\geq 1$, a \emph{linear factor $\calL$ on $\F_p^n$ of complexity $\ell$} is a set of linearly independent vectors  $\calL=\{v_1,\ldots, v_{\ell}\}\subseteq \F_p^n$.   The \emph{atoms of $\calL$} are the sets of the form
$$
L(a_1,\ldots, a_{\ell}):=\{x\in \F_p^n: x^T v_1=a_1,\ldots, x^T v_{\ell}=a_{\ell}\},
$$
for some $a_1,\ldots, a_{\ell}\in \F_p$.  The tuple $a=(a_1,\ldots, a_{\ell})\in \F_p^{\ell}$ is the \emph{label} of the atom $L(a)$.  

By convention, the \emph{linear factor on $\F_p^n$ of complexity $0$} is $\calL=\emptyset$. It has a single atom $\F_p^n$, which we denote by $L(0)$, where $0$ is the unique element in the $0$-dimensional vector space $\F_p^0$.

Given a linear factor on $\F_p^n$, we let $\At(\calL)$ denote the set of atoms of $\calL$.
\end{definition}

Clearly, if two factors $\calL$ and $\calL'$ on $\F_p^n$ satisfy $\Span(\calL)=\Span(\calL')$, then $\At(\calL)=\At(\calL')$.  In the notation of Definition \ref{def:linfac}, $\At(\calL)$ always partitions $\F_p^n$ into $p^{\ell}$ equally sized parts (or atoms). More concretely, the element of $\At(\calL)$ containing $0\in \F_p^n$ is a subgroup of $\F_p^n$, and $\At(\calL)$ consists of the cosets of this subgroup.  On the other hand, given $\ell\geq 0$ and a subspace $H$ of $\F_p^n$ of codimension $\ell$, there exists a linear factor $\calL$ on $\F_p^n$ of complexity $\ell$ so that $\At(\calL)$ is the set of cosets of $H$ (indeed, just take $\calL$ to be any basis of $H^{\perp}$).

Given two linear factors $\calL$ and $\calL'$ on $\F_p^n$, we say that $\calL'$ \emph{refines} $\calL$, denoted $\calL'\preceq \calL$, if the partition $\At(\calL')$ of $\F_p^n$ refines the partition $\At(\calL)$ of $\F_p^n$.

We now define purely quadratic factors.

\begin{definition}[Purely quadratic factor]\label{def:quadfac}
Given an integer $q\geq 1$, a \emph{purely quadratic factor $\calQ$ on $\F_p^n$ of complexity $q$} is a set $\calQ=\{M_1,\ldots, M_q\}$ consisting of $q$ symmetric $n\times n$ matrices over $\F_p$.  The \emph{atoms of $\calQ$} are the sets of the form
$$
Q(b_1,\ldots, b_q):=\{x\in \F_p^n: x^TM_1x=s_1,\ldots, x^TM_qx=s_q\},
$$
for some $b_1,\ldots, b_q\in \F_p$.  The tuple $b\in \F_p^{q}$ is the \emph{label} for the atom $Q(\bbar)$.  

By convention, the \emph{purely quadratic factor of complexity $0$} is $\calQ=\emptyset$. It has a unique atom $\F_p^n$, denoted $Q(0)$, where $0$ is the unique element in the $0$-dimensional vector space $\F_p^0$.

Given a purely quadratic factor $\calQ$ on $\F_p^n$, we let $\At(\calQ)$ denote the set of atoms of $\calQ$. \end{definition}

We now define more general quadratic factors on $\F_p^n$, which  combine linear and purely quadratic factors as follows.

\begin{definition}[Quadratic factor]\label{def:quadfac}
Given integers $\ell,q\geq 0$, a \emph{quadratic factor $\calB$ on $\F_p^n$ of complexity $(\ell,q)$} is a pair $(\calL,\calQ)$ where $\calL$ is a linear factor on $\F_p^n$ of complexity $\ell$ and $\calQ$ is a purely quadratic factor on $\F_p^n$ of complexity $q$.   The \emph{atoms of $\calB$} are the sets of the form
$$
B(a;b):=L(a)\cap Q(b),
$$
for some $a\in \F_p^{\ell}$ and $b\in \F_p^q$.  The tuple $(a,b)\in \F_p^{\ell}\times \F_p^q$ is the \emph{label} for the atom $B(a;b)$.   We let $\At(\calB)$ denote the set of atoms of $\calB$.\end{definition}

Given a quadratic factor $\calB=(\calL,\calQ)$ of complexity $(\ell,q)$, it is clear that we can naturally conflate the set of labels $(a,b)\in \F_p^{\ell}\times \F_p^q$ of the atoms of $\calB$, with the elements of the group $\F_p^{\ell+q}$.  The set of labels, i.e. $\F_p^{\ell+q}$, is then called the \emph{configuration space} associated to the factor $\calB$.  Subsets of the respective configuration spaces  will be denoted by boldface letters.

We will sometimes conflate a linear factor $\calL$ with the quadratic factor $(\calL,\emptyset)$, and similarly, a purely quadratic factor $\calQ$ with the quadratic factor $(\emptyset,\calQ)$.  

Note that when $\calB$ is a quadratic factor of complexity $(\ell,q)$, $\At(\calB)$ forms a partition of $\F_p^n$ with at most $p^{\ell+q}$ parts. However, its parts (or atoms) are not necessarily of the same size. In order to assert that they are at least approximately so, we need to insist that the factor has high rank.

\begin{definition}[Rank of a factor]\label{def:rankfac}
Suppose $\ell\geq 0$ and $q\geq 1$ are itegers, and $\calB=(\calL,\calQ)$ is a quadratic factor on $\F_p^n$ of complexity $(\ell,q)$, where $\calQ=\{M_1,\ldots, M_q\}$.  The \emph{rank of $\calB$} is 
$$
\min\{\rk(\lambda_1M_1+\ldots+\lambda_qM_q): \lambda_1,\ldots, \lambda_q\in \F_p\text{ are not all }0\}.
$$
By convention, the rank of a quadratic factor of the form $(\calL,\emptyset)$ on $\F_p^n$ is $n$.
\end{definition}

The next lemma asserts that when a factor has high rank, all its atoms have approximately the same size (see \cite[Lemma 4.2]{Green.2007}).
\begin{lemma}[Size of an atom]\label{lem:sizeofatoms}
Suppose $\ell,q\geq 0$ are integers, and $\calB=(\calL,\calQ)$ is a quadratic factor on $\mathbb{F}_p^n$ of complexity $(\ell,q)$ and rank $r$.  Then for any $B\in \At(\calB)$,  
\[|B|=(1+O(p^{\ell+q-r/2}))p^{n-\ell-q}.\]
\end{lemma}

We define a \emph{growth function} to be any increasing function $\omega:\mathbb{R}^+\rightarrow \mathbb{R}^+$.  Given two growth functions $\omega, \omega'$, we write $\omega \geq \omega'$ to denote that $\omega(x)\geq \omega'(x)$ for all $x\in \mathbb{R}^+$.  The following immediate corollary of Lemma \ref{lem:sizeofatoms} will be used throughout the paper.

\begin{corollary}\label{cor:sizeofatoms} For any $\mu>0$, there is a polynomial growth function $\rho=\rho(\mu)$ such that if $\calB$ is a quadratic factor on $\mathbb{F}_p^n$ of complexity $(\ell,q)$ and rank at least $\rho(\ell+q)$, then for any $B\in \At(\calB)$, $\big||B|p^{-n}- p^{-(\ell+q)}\big|\leq \mu p^{-(\ell+q)}$.
\end{corollary}

Given two quadratic factors $\calB$ and $\calB'$, we say that $\calB'$ \emph{refines} $\calB$, denoted $\calB'\preceq \calB$, if the partition $\At(\calB')$ of $\F_p^n$ refines the partition $\At(\calB)$ of $\F_p^n$.

\subsection{Reduced colorings}

We begin by defining the reduced coloring of a set in configuration space.

\begin{definition}[Reduced coloring in configuration space associated to $A$ and $\calB$]\label{def:reduced}
Given integers $\ell,q\geq 1$, a set $A\subseteq \mathbb{F}_p^n$, and a quadratic factor $\calB=(\calL,\calQ)$ on $\mathbb{F}_p^n$ of complexity $(\ell,q)$, set 
\begin{align*}
\mathbf{A}_0&=\{a \in \mathbb{F}_p^{\ell+q}: |A\cap B(a)|/|B(a)|\leq \e\},\\
\mathbf{A}_1&=\{a \in \mathbb{F}_p^{\ell+q}: |A\cap B(a)|/|B(a)|\geq 1-\e\},\text{ and }\\
\mathbf{A}_2&=\mathbb{F}_p^{\ell+q}\setminus (\mathbf{A}_1\cup \mathbf{A}_0).
\end{align*}
The \emph{reduced coloring associated to $A$ and $\calB$} is the $3$-coloring
\[\Red_{\calB,\e}(\F_p^n,A):=(\mathbb{F}_p^{\ell+q};\mathbf{A}_0,\mathbf{A}_1,\mathbf{A}_2).\]
\end{definition}

The sets $\mathbf{A}_0,\mathbf{A}_1,\mathbf{A}_2$ defined above depend on $\calB$, $A$ and $\e$, but we omit this dependence from the notation as it will always be clear from the context.  

Note that in the case where $A$ has bounded $\VC_2$-dimension, and $\calB$ arises from Theorem \ref{thm:vc2}, there will be very few atoms of $\calB$ of intermediate density. In this case, the set $\mathbf{A}_2$ will be small, and the reduced coloring $\Red_{\calB,\e}(\F_p^n,A)$ will encode useful structural information about the set $A$.  

Given $\Red_{\calB,\e}(\F_p^n,A)=(\F_p^{\ell+q}; \mathbf{A}_0, \mathbf{A}_1,\mathbf{A}_2)$, we will consider the VC-dimension of $\mathbf{A}_1|\mathbf{A}_0$, relative to a  special subspace of $\F_p^{\ell+q}$ defined using the linear part of the factor $\calB$. 

\begin{definition}[$\bH_{\calB}$ and $\calL_{\calB}$]\label{def:hb}
Let $\ell,q,n\geq 1$ be integers, and suppose that $\calB=(\calL,\calQ)$ is a quadratic factor on $\F_p^n$ of complexity $(\ell,q)$. Define $\calL_{\calB}=\{e_1,\ldots, e_{\ell}\}$, where $e_i$ is the $i$th standard basis vector in $\F_p^{\ell+q}$, and define 
$$
\bH_{\calB}=\{0 b\in \F_p^{\ell+q}: b\in \F_p^q\}=\{d\in \F_p^{\ell+q}: d_1=\ldots=d_{\ell}=0\}.
$$ 
\end{definition}

Observe that in the notation of Definition \ref{def:hb}, the set $\bH_{\calB}$ is a subspace in $\F_p^{\ell+q}$ of codimension $\ell$, and is also the zero atom of the factor $\calL_{\calB}$ on $\F_p^{\ell+q}$. 

A key ingredient in the proof of our main result is Theorem \ref{thm:extract3} below, which says that if $ \mathbf{A}_1|\mathbf{A}_0$ has VC-dimension at least $k$ relative to $\bH_{\calB}$, then $A$ has  $k$-$\IP_2$.  The ``relative to $\bH_{\calB}$"  assumption is key, as it forces the copy of $k$-$\IP$ to use the purely quadratic structure in a nontrivial way (see Section 4 of \cite{Terry.2021a} for further discussion on this).  

\begin{theorem}[Sufficient condition for $\IP_2$]\label{thm:extract3}
For all primes $p>2$, integers $k\geq 1$, and reals $\e\in (0,1)$, there exists a polynomial growth function $\omega=\omega(p,k,\e)$ such that for all integers $\ell, q\geq 1$ and all sufficiently large $n$, the following holds.  

Let $A\subseteq \mathbb{F}_p^n$, and let $\calB=(\calL,\calQ)$ be a quadratic factor on $\F_p^{\ell+q}$ of complexity $(\ell,q)$ and rank at least $\omega(\ell+q)$.  Suppose that $\mathbf{A}_1|\mathbf{A}_0$ has $\VC$-dimension at least $k$ relative to $\bH_{\calB}$, where $\Red_{\calB,\e}(\F_p^n,A)=(\F_p^{\ell+q};\mathbf{A}_0,\mathbf{A}_1, \mathbf{A}_2)$.  Then $A$ has $k$-$\IP_2$.
\end{theorem}

Theorem \ref{thm:extract3} is proved as Proposition 4.20 in \cite{Terry.2021a}, where it serves as a companion to an analogous result relating $\FOP_2$, a higher order generalisation of stability introduced in that paper, to copies of half-graphs with one side in $\bH_{\calB}$.  

To conclude this section, we prove a corollary of our approximate local version of the Alon-Fox-Zhao theorem (Theorem \ref{thm:errorstructure}). For this, it will be helpful to have the following terminology at our disposal.

\begin{definition}[Almost homogeneity]\label{def:homog}
Given a partition $\calP$ of a group $G$ and a set $A\subseteq G$, we say that $\calP$ is \emph{almost $\mu$-homogeneous with respect to $A$ } if 
$$
\Big|\bigcup_{P\in \Gamma}P\Big|\geq (1-\mu)|G|,
$$
where $\Gamma=\{P\in \calP: \text{ either $|A\cap P|\geq (1-\mu)|P|$ or $|A\cap P|\leq \mu|P|$}\}$.  

A linear or quadratic factor $\calB$ on $\F_p^n$ is \emph{almost $\mu$-homogeneous with respect to $A$} if the partition $\At(\calB)$ is almost $\mu$-homogeneous with respect to $A$. 
\end{definition}

\begin{corollary}\label{cor:extract1}
For all primes $p>2$ and integers $k\geq 0$, there exists $\delta_0\in (0,1)$ so that for all $0<\delta<\delta_0$ and $0<\e\leq (\delta/120)^{k+1}$, there is a growth function $\omega=\omega(p,k,\delta,\e)$ such that for all integers $\ell, q\geq 1$, the following holds for all sufficiently large $n$.  

Let $A\subseteq \mathbb{F}_p^n$ have $\VC_2$-dimension at most $k$, and suppose $\calB=(\calL,\calQ)$ is a quadratic factor on $\F_p^n$ which has complexity $(\ell,q)$, which has  rank at least $\omega(\ell+q)$, and which is almost $\e$-homogeneous with respect to $A$. Let $\Red_{\calB,\e}(\F_p^n,A)=(\F_p^{\ell+q};\mathbf{A}_0,\mathbf{A}_1, \mathbf{A}_2)$.  

Then there is a subspace $\bH'\leqslant \bH_{\calB}$ of index $\delta^{-k-o_{k,p}(1)}$ in $\bH_{\calB}$ and unions $\bY_0,\bY_1$ of cosets of $\bH'$ in $\F_p^{\ell+q}$ so that for each $u\in \{0,1\}$, 
$$
|\mathbf{A}_u\Delta \bY_u|\leq 4\delta^{1-o_{k,p}(1)}|\F_p^{\ell+q}|,
$$
 where $o_{k,p}(1)$ tends to $0$ as $\e$ tends to $0$ at a rate depending on $k$ and $p$.  
 
 Moreover, $o_{k,p}(1)$ can be taken to have the form $\psi(p,k,\delta)$ from Theorem \ref{thm:errorstructure}.
\end{corollary}

\begin{proof}
Let $\delta_0=\delta_0(p,k)$ from Theorem \ref{thm:errorstructure}. Fix $0<\delta<\delta_0$, and let $\psi(p,k,\delta)$ be as in Theorem \ref{thm:errorstructure}.  Fix $0<\e\leq (\delta/120)^{k+1}$,  and let $\omega=\omega(p,k,\e)$ be as in Theorem \ref{thm:extract3}.  Fix integers $\ell,q\geq 1$, and let $n$ be sufficiently large.

Fix $A\subseteq \mathbb{F}_p^n$ with $\VC_2$-dimension at most $k$. Suppose $\calB=(\calL,\calQ)$ is a quadratic factor on $\F_p^n$ of complexity $(\ell,q)$ and rank at least $\omega(\ell+q)$ that is almost $\e$-homogeneous with respect to $A$.  By Theorem \ref{thm:extract3}, $\mathbf{A}_1|\mathbf{A}_0$ has VC-dimension at most $k$ relative to $\bH_{\calB}$, where $\Red_{\calB,\e}(\F_p^n,A)=(\F_p^{\ell+q};\mathbf{A}_0,\mathbf{A}_1, \mathbf{A}_2)$.  By Theorem \ref{thm:errorstructure}, there is a subspace $\bH'\leqslant \bH_{\calB}$ of index at most $\delta^{- k-\psi(p,k,\delta)}$  in $\bH_{\calB}$ so that the desired conclusion holds.  
\end{proof}

\section{Proofs of Theorem \ref{thm:vc2main} and Corollary \ref{cor:betterbound}}\label{sec:betterproof}
Corollary \ref{cor:vc2} tells us that given a set $A\subseteq \mathbb{F}_p^n$ of $\VC_2$-dimension at most $k$, we can find a quadratic factor $\calB$ of high rank so that on most atoms of $\calB$, the density of $A$ is either close to $0$ or close to $1$, meaning that $\calB$ is almost homogeneous with respect to $A$ (in the sense of Definition \ref{def:homog}). By Corollary \ref{cor:extract1} at the end of the preceding section, we can now find a bounded index subspace of $\bH_{\calB}$, whose cosets can be used to approximate the sets $\mathbf{A}_0$ and $\mathbf{A}_1$ from  the reduced coloring $\Red_{\calB,\e}(\F_p^n,A)$.  The last remaining step is to show that the partition of configuration space given by this subspace can be pulled back to a quadratic factor on the original group $\F_p^n$, which can then be used to approximate the original set $A$. The main idea for doing so was developed in our earlier paper \cite{Terry.2021d}, where we defined the following (see \cite[Definition 5.13]{Terry.2021d}).

\begin{definition}[Pullback of a linear factor on configuration space]\label{def:pullback}
Let $p>2$ be a prime, and let $\ell,q,n\geq 1$ be integers.  Suppose $\calB=(\calL,\calQ)$ is a quadratic factor on $\F_p^n$ of complexity $(\ell,q)$, and $\calR$ is a linear factor on $\F_p^{\ell+q}$.  For each $\bR\in \At(\calR)$, set $X_{\bR}=\bigcup_{a\in \bR}B(a)$, where each $B(a)$ is an atom of $(\calL,\calQ)$.  Define the \emph{pullback $\calX_{\calR}$ of $\calR$} to be
 $$
 \calX_{\calR}=\{X_{\bR}: \bR\in \At(\calR)\}.
 $$
\end{definition}
Note that $X_{\bR}\subseteq \F_p^n$ for each $\bR\in \At(\calR)$.  Clearly $\calX_{\calR}$ partitions $\F_p^n$, since $\At(\calR)$ partitions $\At(\calB)$ and $\At(\calB)$ partitions $\F_p^n$.  Thus Definition \ref{def:pullback} gives us a way of building a partition $\calX_{\calR}$ of $\F_p^n$ from the atoms of the linear factor $\calR$ on the configuration space associated to $\calB$, i.e. $\F_p^{\ell+q}$.  

Our next lemma gives a sufficient condition for when a $A\subseteq\F_p^n$ can be well approximated by a union of sets from a partition of the form $\calX_{\calR}$ from Definition \ref{def:pullback}

\begin{lemma}[Pullbacks preserve approximations]\label{lem:pullpres}
Let $p>2$ be a prime.  For all reals $0<\mu<1/3$ and $0<\e<\mu/4$, there is a polynomial growth function $\rho=\rho(p,\mu)$ such that for all integers $\ell,q\geq 1$ and all sufficiently large $n$, the following holds.

Let $A\subseteq \F_p^n$, and suppose $\calB=(\calL,\calQ)$ is a quadratic factor on $\F_p^n$ of complexity $(\ell,q)$ and rank at least $\rho(\ell+q)$ that is almost $\e$-homogeneous with respect to $A$. Let $\Red_{\calB,\e}(\F_p^n,A)=(\F_p^{\ell+q};\mathbf{A}_0,\mathbf{A}_1,\mathbf{A}_2)$.  

Suppose there exists a linear factor $\calR$ on $\F_p^{\ell+q}$ along with unions $\bY_0,\bY_1$ of atoms of $\calR$ with the property that for each $u\in \{0,1\}$, $|\mathbf{A}_u\Delta \mathbf{Y}_u|\leq \mu |\F_p^{\ell+q}|$.  Then there is a union $Y$ of sets from $\calX_R$ such that $|A\Delta Y|\leq 4\mu$.  
\end{lemma}
\begin{proof}
Fix $p>2$ a prime and reals $0<\mu<1/3$ and $0<\e<\mu/4$.  Let $\rho=\rho(p,\mu/2)$ be from Corollary \ref{cor:sizeofatoms}. Fix integers $\ell,q\geq 1$, and let $n$ be sufficiently large. 

Suppose $A\subseteq \F_p^n$ and $\calB=(\calL,\calQ)$ is a quadratic factor on $\F_p^n$ which has complexity $(\ell,q)$, which has rank at least $\rho(\ell+q)$, and which is almost $\e$-homogeneous with respect to $A$. Let $\Red_{\calB,\e}(\F_p^n,A)=(\F_p^{\ell+q};\mathbf{A}_0,\mathbf{A}_1,\mathbf{A}_2)$.

Let $\calR$ be a linear factor on $\F_p^{\ell+q}$ with the property that for some   unions $\bY_0,\bY_1$ of atoms of $\calR$, $|\mathbf{A}_u\Delta \bY_u|\leq \mu |\F_p^{\ell+q}|$ holds for each $u\in \{0,1\}$.  Let $r\geq 0$ be the complexity of $\calR$, and note this implies that every $\bR\in \At(\calR)$ has size $p^{\ell+q-r}$.  Since $\bY_0,\bY_1$ are unions of atoms of $\calR$, there exist $Y_0,Y_1\subseteq \At(\calR)$ so that $\bY_0=\bigcup_{\mathbf{R}\in Y_0}\mathbf{R}$ and $\bY_1=\bigcup_{\mathbf{R}\in Y_1}\mathbf{R}$.  Consider the subset $Y$ of $\F_p^n$ defined by
$$
Y:=\bigcup_{\mathbf{R}\in Y_1}X_\mathbf{R},
$$
i.e. by taking the union of those sets from $\calX_{\calR}$ indexed by elements of $Y_1$.

We show that $Y$ is our desired approximation for the original set $A$.  We will use the fact that Corollary \ref{cor:sizeofatoms} and our choice of rank function imply 
\begin{align}\label{bdbound}
\text{ for every $B(d)\in \At (\calB)$, }|B(d)|=(1\pm \mu/2)p^{n-\ell-q}.
\end{align}

Observe
\begin{align*}
|A\setminus Y|=\sum_{\bR\in \At(\calR)\setminus Y_1}|A\cap X_\bR|&=\sum_{\bR\in \At(\calR)\setminus Y_1}\Big(\sum_{d\in \bR}|A\cap B(d)|\Big)\\
&\leq \sum_{\bR\in \At(\calR)\setminus Y_1}\Big(\sum_{d\in \bR\cap \mathbf{A}_0}|A\cap B(d)|+\sum_{d\in \bR\cap  \mathbf{A}_2}|B(d)|+\sum_{d\in \bR\cap \mathbf{A}_1}|B(d)|\Big)\\
&\leq \sum_{d\in \mathbf{A}_0}|A\cap B(d)|+\sum_{d\in  \mathbf{A}_2}|B(d)|+\sum_{d\in (\mathbf{A}_1\setminus \bY_1)}|B(d)|\\
&\leq \sum_{d\in \mathbf{A}_0}\e|B(d)|+\sum_{d\in  \mathbf{A}_2}|B(d)|+\sum_{d\in (\mathbf{A}_1\setminus \bY_1)}|B(d)|,
\end{align*}
where the last inequality uses the definition of $\mathbf{A}_0$.  Using (\ref{bdbound}), $|A\setminus Y|$ is at most
\begin{align*}
&\e|\mathbf{A}_0|(1+\mu/2)p^{n-\ell-q}+|\mathbf{A}_2|(1+\mu/2)p^{n-\ell-q}+|\mathbf{A}_1\setminus \bY_1|(1+\mu/2)p^{n-\ell-q}\\
&\leq \e (1+\mu/2)p^n+\e (1+\mu/2)p^n+\mu (1+\mu/2)p^n\\
&\leq 2\mu p^n,
\end{align*}
where the second inequality uses the fact that $\calB$ is almost $\e$-homogeneous with respect to $A$ (which implies that $|\mathbf{A}_2|\leq \e p^{\ell+q}$) and our assumption that $|\mathbf{A}_1\Delta \bY_1|\leq \mu p^{\ell+q}$, and the final inequality holds because $\e<\mu/2$ and $\mu$ is sufficiently small. We have thus shown that $|A\setminus Y|\leq 2\mu |G|$.  Similarly, we have 
\begin{align*}
|Y\setminus A|=\sum_{\bR\in Y_1}|X_{\mathbf{R}}\setminus A|&\leq \sum_{\bR\in Y_1}\Big(\sum_{d\in \bR\cap \mathbf{A}_0}|B(d)|+\sum_{d\in \bR\cap \mathbf{A}_2}|B(d)|+\sum_{d\in \bR\cap \mathbf{A}_1}| B(d)\setminus A|\Big)\\
&\leq \sum_{d\in \bY_1\setminus \mathbf{A}_1}|B(d)|+ \sum_{d\in  \mathbf{A}_1}|B(d)\setminus A|\\
&\leq \sum_{d\in \bY_1\setminus \mathbf{A}_1}|B(d)|+ \sum_{d\in  \mathbf{A}_1}\e|B(d)|,
\end{align*}
where the last inequality uses the definition of $\mathbf{A}_1$.  Using (\ref{bdbound}), $|Y\setminus A|$ is therefore at most 
\begin{align*}
|\bY_1\setminus \mathbf{A}_1|(1+\mu/2)p^{n-\ell-q}+\e|\mathbf{A}_1|(1+\mu/2)p^{n-\ell-q}\leq \mu (1+\mu/2)p^n+\e (1+\mu/2)p^n\leq 2\mu p^n,
\end{align*}
where the first inequality uses the assumption that $| \mathbf{A}_1\Delta \bY_1|\leq \mu |\F_p^{\ell+q}|$, and the second holds because $\e<\mu/4$ and $\mu$ is sufficiently small.  We have therefore shown that $|Y\setminus A|\leq 2\mu |G|$.  Combining, we have that $|A\Delta Y|\leq 4\mu |G|$, as desired.
\end{proof}

Given partitions of the form $\calX_{\calR}$ as in Definition \ref{def:pullback} and $\calL_{\calB}$ is as in Definition \ref{def:hb}, if we moreover know that $\calR\preceq\calL_{\calB}$,  then $\calX_{\calR}$ is equal to the set of atoms of a quadratic factor on $\F_p^n$.   In particular, we proved the following in \cite[Lemma 4.24]{Terry.2021a}.   

 \begin{lemma}[Extracting a quadratic factor \cite{Terry.2021a}]\label{lem:extractfactor}
Let $\ell,q,n\geq 1$ be integers and let $\rho>0$ be a real.  Suppose  $\calB=(\calL,\calQ)$ is a quadratic factor on $\F_p^n$ of complexity $(\ell,q)$ and rank $\rho$, and $\calR\preceq\calL_{\calB}$ is a linear factor on $\F_p^{\ell+q}$ of complexity $\ell+r$.  Then there is a quadratic factor $\tilde{\calB}=(\calL,\tilde{\calQ})$ on $\F_p^n$ of complexity $(\ell, r)$ and rank at least $\rho$ such that $\At(\tilde{\calB})=\calX_{\calR}$.
\end{lemma}

We can combine this with Corollary \ref{cor:extract1} to show that given a high rank quadratic factor which is homogeneous with respect to a set $A$ of bounded $\VC_2$-dimension, there exists another factor $\calB'$ on $\F_p^n$ which has the same linear component, which has a possibly new quadratic component of very small size, and such that $A$ is approximately a union of atoms of $\calB'$.

\begin{theorem}\label{thm:extract1}
For all primes $p>2$ and integers $k\geq 0$, there exists $\delta_0=\delta_0(k)>0$ so that for all $0<\delta<\delta_0$ and all $0<\e\leq (\delta/120)^{k+2}$, there is a polynomial growth function $\omega=\omega(p,k,\delta,\e)$ such that for all integers $\ell, q\geq 1$ the following holds for all sufficiently large $n$.

Let $A\subseteq \mathbb{F}_p^n$ be a subset of $\VC_2$-dimension at most $k$, and suppose that $\calB=(\calL,\calQ)$ is a quadratic factor on $\F_p^n$ of complexity $(\ell,q)$ and rank at least $\omega(\ell+q)$ that is almost $\e $-homogeneous with respect to $A$.   Then there exist an integer $q'\geq 0$ and a quadratic factor $\calB'=(\calL,\calQ')$ on $\F_p^n$ of complexity $(\ell,q')$ and rank at least $\omega(\ell+q')$ such that  
\begin{enumerate}[label=\normalfont(\roman*)]
\item $q'\leq \log_p( \delta^{-k-o_{k,p}(1)})$;
\item there exists a union $Y$ satisfying $|A\Delta Y|\leq 16\delta^{1-o_{k,p}(1)}|G|$;
\end{enumerate}
where $o_{k,p}(1)$ tends to $0$ as $\delta$ tends to $0$ at a rate depending on $k$ and $p$.  Moreover, $o_{k,p}(1)$ can be taken to have the form $\psi(k,p,\delta)$ from Theorem \ref{thm:errorstructure}s and \ref{thm:extract3}.
\end{theorem}
\begin{proof}
Fix a prime $p>2$ and an integer $k\geq 0$, and let $\delta_0=\delta_0(p,k)$ be as in Corollary \ref{cor:extract1}. Fix $0<\delta<\delta_0$ and $0<\e\leq (\delta/120)^{k+1}$. Let $\psi(p,k,\delta)$ be as in Corollary \ref{cor:extract1}.  Let $\omega_0=\omega_0(p,k,\delta)$ be a growth function as in Corollary \ref{cor:extract1}, and let $\omega_1=\omega_1(p,k,\e)$ be a growth function as in Theorem \ref{thm:extract3}.  Now let $\omega$ be the pointwise maximum of $\omega_0$ and $\omega_1$. Fix integers $\ell,q\geq 1$   and $n$ sufficiently large.  

Let $A\subseteq \mathbb{F}_p^n$ have $\VC_2$-dimension at most $k$. Suppose that $\calB=(\calL,\calQ)$ is a quadratic factor of complexity $(\ell,q)$ and rank at least $\omega(\ell+q)$ which is almost $\e$-homogeneous with respect to $A$.   By Theorem \ref{thm:extract3}, $\mathbf{A}_1|\mathbf{A}_0$ has $\VC$-dimension at most $k$ relative to $\bH_{\calB}$, where $\Red_{\calB,\e}(\F_p^n,A)=(\F_p^{\ell+q};\mathbf{A}_0,\mathbf{A}_1, \mathbf{A}_2)$.

 By Corollary \ref{cor:extract1}, there exist an integer $r\leq \log_p(\delta^{-k-\psi(p,k,\delta)})$, a subspace $\bH'\leqslant \bH_{\calB}$ of codimension $r$ in $\bH_{\calB}$, and unions $\bY_0,\bY_1$ of cosets of $\bH'$ in $\F_p^{\ell+q}$ so that for each $u\in \{0,1\}$, $|\mathbf{A}_u\Delta \bY_u|\leq 4\delta^{1-\psi(p,k,\delta)} |\F_p^{\ell+q}|$.  Let $\calR$ be a linear factor on $\F_p^{\ell+q}$ such that the atoms of $\calR$ are exactly the cosets of $\bH'$.  Since $\bH'\leqslant \bH_{\calB}$, $\calR\preceq \calL_{\calB}$, and $\calR$ must have complexity $r+\ell$. Applying Lemma \ref{lem:pullpres} with $\mu=4\delta^{1-\psi(p,k,\delta)}$, we obtain a union $Y$ of sets from $\calX_{\calR}$, with the property that $|A\Delta Y|\leq 16\delta^{1-\psi(p,k,\delta)}|G|$.  By Lemma \ref{lem:extractfactor},  there is a quadratic factor $\tilde{\calB}=(\calL,\tilde{\calQ})$ on $\F_p^n$ with $\At(\tilde{\calB})=\calX_{\calR}$, whose complexity is $(\ell, r)$, and whose rank is at least $\omega(\ell+q)$.  Note $\calR$ can have complexity no larger than $\ell+q$, meaning $\ell+r\leq \ell+q$. Consequently, we have $\rk(\tilde{\calB})\geq \omega(\ell+q)\geq \omega(\ell+r)$. This finishes the proof.
\end{proof}

We now put everything together to prove our main structure theorem, stated as Theorem \ref{thm:vc2main} in the introduction.

\begin{theorem}\label{thm:vc2alt}
For all primes $p>2$ and integers $k\geq 0$, there is $\mu_0\in (0,1)$ so that for all $0<\mu<\mu_0$, there exists a growth function $\omega_0=\omega_0(p,k,\mu)$ such that for all $\omega\geq\omega_0$ there exists a constant $C=C(k,\mu,\omega)>0$ such that the following holds for all sufficiently large $n$.

Suppose that $A\subseteq \F_p^n$ satisfies $\VC_2(A)\leq k$.  Then there exist integers $\ell,q'\geq 1$ and a quadratic factor $\calB'$ on $\F_p^n$ of complexity $(\ell,q')$ such that
\begin{enumerate}[label=\normalfont(\roman*)]
\item $\ell\leq C$ and $q'\leq \log_p(\mu^{-k-o_{k,p}(1)})$;
\item $\calB'$ has rank at least $\omega(\tilde{\ell}+\tilde{q})$;
\item there exists a union $Y$ of atoms $\calB'$ satisfying $|A\Delta Y|\leq \mu |G|$;
\end{enumerate}
  where $o_{k,p}(1)$ denotes a function $\psi'(k,p,x)$ which tends to $0$ as $x$ tends to $0$, at a rate depending on $k$ and $p$.  
\end{theorem}

\begin{proof}
Fix a prime $p>2$ and an integer $k\geq 0$. Let $\delta_0=\delta_0(k)$ and $\psi(p,k,x)$ be as in Theorem \ref{thm:extract1}. Recall $\psi(p,k,x)$ has the form $\psi(p,k,x)=\sigma_{p,k}\log(x^{-1})^{-1/5}$, where $\sigma_{p,k}$  is a constant depending on $p$ and $k$.  Given $x\in (0,1/2)$, define 
$$
\Theta(x)=16x^{1-\psi(p,k,x)}.
$$
Clearly $\Theta(x)\rightarrow 0$ as $x\rightarrow 0$ (since $\psi(p,k,x)$ does).  It is not difficult to show $\Theta$ has an inverse, considered as a function from $(0,1/2)$ into $(0,\Theta(1/2))$ (since it is differentiable and increasing).  Note that $\Theta^{-1}(x)\rightarrow 0$ as $x\rightarrow 0$. Choose $\mu_0\in (0,1/2)$ sufficiently small so that $\mu_0<\Theta(\delta_0)$, and so that $x<\Theta^{-1}(\mu)$ implies $\psi(p,k,x)\leq 1/2$.

 Fix $0<\mu<\mu_0$, and let $\delta=\Theta^{-1}(\mu)$. Since $\Theta$ is increasing, $\delta<\delta_0$.  Set $\e=(\delta/120)^{k+1}$.  Let $\rho_0=\rho_0(k,\delta,\e)$ be from Theorem  \ref{thm:extract1}, let $\rho_1=\rho_1(p,k,\e)$ be from Theorem \ref{thm:vc2}, and let $\omega_0$ be the pointwise maximum of $\rho_0$ and $\rho_1$.  Fix a growth function $\omega \geq \omega_0$, let $C=C(p,k,\e,\omega)$  be from Theorem \ref{thm:vc2}, and let $n$ be sufficiently large.  

Suppose $A\subseteq \F_p^n$ satisfies $\VC_2(A)\leq k$. Apply Theorem \ref{thm:vc2} to find integers $1\leq \ell,q\leq C$ and a quadratic factor $\calB=(\calL,\calQ)$ on $\F_p^n$ of complexity $(\ell,q)$ and rank at least $\omega(\ell+q)$ that is almost $\e$-homogeneous with respect to $A$.  By Theorem \ref{thm:extract1}, there exist an integer $0\leq q'\leq \log_p( \delta^{-k-\psi(p,k,\delta)})$, a quadratic factor $\calB'=(\calL,\calQ')$ on $\F_p^n$ of complexity $(\ell,q')$ and rank at least $\omega(\ell+q')$, and a union $Y$ of atoms of $\calB'$ such that 
$$
|A\Delta Y|\leq 16\delta^{1-\psi(p,k,\delta)}|G|=\mu|G|.
$$
We just have left to prove the desired bound on $q'$. Since we already know $q'\leq \log_p( \delta^{-k-\psi(p,k,\delta)})$, it suffices to show $ \delta^{-k-\psi(p,k,\delta)}\leq \mu^{-k-o_{k,p}(1)}$, where $o_{k,p}(1)\rightarrow 0$ as $\mu\rightarrow 0$, at a rate depending on $k$ and $p$.  To this end, observe that since $\delta=(\mu/16)^{(1-\psi(p,k,\delta))^{-1}}$ and $\psi(p,k,\delta)\leq 1/2$,
$$
\delta^{-k-\psi(p,k,\delta)}=\Big(\frac{\mu}{16}\Big)^{\frac{-k-\psi(p,k,\delta)}{1-\psi(p,k,\delta)}}=\Big(\frac{\mu}{16}\Big)^{-k-\frac{(k+1)\psi(p,k,\delta)}{1-\psi(p,k,\delta)}}\leq\Big(\frac{\mu}{16}\Big)^{-k-2(k+1)\psi(p,k,\delta) }=\mu^{-k-\psi'(p,k,\mu)},
$$
where, recalling $\delta=\Theta^{-1}(\mu)$, we have 
$$
\psi'(p,k,\mu)=2(k+1)\psi(p,k,\Theta^{-1}(\mu)) +(k+2(k+1)\psi(p,k,\Theta^{-1}(\mu)) )\log_{\mu}(16). 
$$
  Since $\Theta(\mu^{-1})\rightarrow 0$ as $\mu \rightarrow 0$, $\psi(p,k,\Theta(\mu^{-1}))\rightarrow 0$ as $\mu\rightarrow 0$, and consequently, $\psi'(p,k,\mu)\rightarrow 0$ as $\mu\rightarrow 0$. This finishes the proof.
\end{proof}

In order to deduce Corollary \ref{cor:betterbound}, we use the following standard averaging fact.

\begin{lemma}\label{lem:avcor}
Let $\e,a,b\in (0,1)$ be such that $\e=ab$. Let $X$ be a finite set, and let $A\subseteq X$ satisfy $|A|\leq \e|X|$.  Given any partition $\calP$ of $X$, define $\Sigma=\{P\in \calP: |A\cap P|> a|P|\}$. Then $|\bigcup_{P\in \Sigma}P|\leq b|X|$. 
\end{lemma}

We can now deduce our main result, Corollary \ref{cor:betterbound}. 

\vspace{3mm}

\begin{proofof}{Corollary \ref{cor:betterbound}}
Fix a prime $p>2$ and an integer $k\geq 0$, and let $\mu_0=\mu_0(p,k)$ be from Theorem \ref{thm:vc2alt}.  Given $x\in (0,1/2)$, define $\varphi(x)=x^{1+\log_x((1-x)^{-1})}$.  It is not difficult to check $\varphi$ has an inverse, considered as a function from $(0,1/2)$ into $(0,\varphi(1/2))$ (since it is differentiable and increasing).  Note $\varphi^{-1}(x)\rightarrow 0$ as $x\rightarrow 0$.  Let $\e_0=\varphi(\mu_0)$ and fix $0<\e<\e_0$.  Let $\e_1=(\e/2)^8$, and set $\mu=\varphi^{-1}(\e_1)$. Note $0<\mu<\mu_0$, since $\e_1<\e$ and $\varphi$ is increasing.

Let $\omega_1=\omega_1(k,\mu^{2})$ be as in Theorem \ref{thm:vc2alt},  let $\omega_2=\omega_2(\mu)$ be as in Corollary \ref{cor:sizeofatoms}, and let $\omega_0$ be the pointwise maximum of $\omega_1$ and $\omega_2$.  Fix a growth function $\omega\geq \omega_0$, and let  $C_0=C_0(k,\mu^{2},\omega)$ be as in Theorem \ref{thm:vc2alt}.   Let $\psi'(p,k,x)$ be the $o_{k,p}(1)$ function from Theorem \ref{thm:vc2alt}.  Note that since $\varphi(\mu)=\e_1$, $\mu<\e_1$.

Assume now that $n$ is sufficiently large and $A\subseteq \F_p^n$ has $\VC_2$-dimension at most $k$.  By Theorem \ref{thm:vc2alt}, there exist integers $\ell,q$ satisfying $1\leq \ell \leq C_0$ and $0\leq q\leq  \log_p(\mu^{-2k-2\psi'(p,k,\mu^2)})$, a quadratic factor $\calB=(\calL,\calQ)$ on $\F_p^n$ of complexity $(\ell,q)$ and rank at least $\omega(\ell,q)$, and a union $Y$ of atoms of $\calB$ so that $|A\Delta Y|\leq \mu^2
|G|$.  Define
$$
\Xi=\{B\in \At(\calB): |(A\Delta Y)\cap B|>\mu|B|\}.
$$
  By Lemma \ref{lem:avcor}, $|\bigcup_{B\in \Xi}B|\leq \mu |G|$.  Combining with Corollary \ref{cor:sizeofatoms}, we have  
$$
|\Xi|(1-\mu^2)p^{n-\ell-q}\leq \Big|\bigcup_{B\in \Xi}B\Big|\leq \mu |G|,
$$
and consequently, $|\Xi|\leq \mu(1-\mu)^{-1}p^{\ell+q}=\varphi(\mu)p^{\ell+q}=\e_1 p^{\ell+q}\leq \e p^{\ell+q}$. Suppose now $B\in \At(\calB)\setminus \Xi$.  By definition of $\Xi$, $|(A\Delta Y)\cap B|\leq \mu |B|<\e_1|B|$. If $B\in Y$, then 
$$
|(A\Delta Y)\cap B|=|B\setminus A|\leq \mu |B|<\e_1|B|.
$$
 On the other hand, if $B\notin Y$, then $|(A\Delta Y)\cap B|=|B\cap A|\leq \mu |B|<\e_1|B|$.  Now let $\Gamma$ contain all tuples $d\in \F_p^{3\ell+6q}$ so that $\Sigma_{\calB}(d)\notin \Xi$.  For each $d\in \Gamma$, we have shown above that $\alpha_{B(\sigma_{\calB}(d))}\in [0,\e_1)\cup (1-\e_1,1]$, so by Fact \ref{fct:unihom2} and definition of $\e_1$, we have that for all $d\in \Gamma$, $\|1_A-\alpha_{B(\Sigma_{\calB}(d))}\|_{U^3(d)}<2\e_1^{1/8}=\e$.  Since $|\Xi|\leq \e p^{\ell+q}$, we have $|\Gamma|\geq (1-\e)p^{3\ell+6q}$.  
 
We now consider the bound on the parameter $q$. Recalling that $\phi(\mu)=\e_1$ and $\e_1=2\e^{1/8}$, we have 
\begin{align*}
q\leq \log_p(\mu^{-2k-2\psi'(p,k,\mu^2)} )&=\log_p((\e_1(1-\mu))^{-2k-2\psi'(p,k,\mu^2)})\\
&\leq \log_p(\e_1^{-2k-2\psi'(p,k,\mu^2)})\\
&=\log_p(\e^{-16k-16\psi'(p,k,\varphi^{-1}(2\e^{1/8})^2)-(2k+2\psi'(p,k,\varphi^{-1}(2\e^{1/8})^2))\log_\e(2)})\\
&=\log_p(\e^{-16k-o_{k,p}(1)}), 
\end{align*}
where $o_{k,p}(1)$ tends to $0$ as $\e$ tends to $0$ (because $\varphi^{-1}(2\e^{1/8})^2\rightarrow 0$ as $\e\rightarrow 0$, and thus $\psi'(p,k,\varphi^{-1}(2\e^{1/8})^2)\rightarrow 0$ as $\e \rightarrow 0$).
\end{proofof}

\begin{remark}\label{rem:strength}
As noted in the introduction, given a quadratic factor that is sufficiently uniform with respect to a set $A$ of bounded $\VC_2$-dimension (in the sense of Theorem \ref{thm:vc2}), the density of $A$ is close to $0$ or $1$ on all uniform atoms (note that all parameter relationships in Fact \ref{fct:unihom2} are polynomial). It is therefore possible to prove a version of Theorem \ref{thm:extract1} in which one assumes that the factor is merely sufficiently uniform (rather than homogeneous) with respect to $A$. Thus, if one managed to improve the bound on the linear part of the factor in Theorem \ref{thm:vc2}, then one could use the above argument to improve the bound on the quadratic complexity while keeping the same linear component.
\end{remark}

\section{A lower bound}\label{sec:lowerbound}
Our goal in this section is to prove Theorem \ref{thm:lowerbound} by giving an example of a subset of $\F_p^n$ for which Corollary \ref{cor:vc2} cannot hold. Specifically, the example has the property that any quadratic factor that provides a $U^3$ regularity decomposition with uniformity parameter $\e$ must have quadratic complexity at least a power of $\e^{-1}$. It is a separate open problem (which we do not address here) to produce an example showing that either $\ell$ or $q$ must match the tower type upper bounds arising from the proof of Corollary \ref{cor:vc2}).

 For the analysis of our example, we will require the following fact, proved in \cite{Terry.2021a} (see Proposition A.22 there).  

\begin{proposition}\label{prop:preservedensity}
 For all primes $p>2$ and $\e\in (0,1)$, there is a polynomial growth function $\rho=\rho(p,\e)$ such that the following holds.  Suppose $\ell,q\geq 1$ are integers, $A\subseteq \F_p^n$, $\calB=(\calL,\calQ)$ is a quadratic factor on $\F_p^n$ of complexity $(\ell, q)$ and rank at least $\rho(\ell+q)$, and 
 $$
 d=(a_1,a_2,a_3,b_{12},b_{13},b_{23})\in \F_p^{\ell+q}\times \F_p^{\ell+q}\times\F_p^{\ell+q}\times\F_p^{q}\times\F_p^{q}\times\F_p^{q}
 $$
  is such that $\|1_A-\alpha_{B(\Sigma(d))}\|_{U^3(d)}<\e$.  
 
 Suppose $\calB'=(\calL',\calQ')$ is a quadratic factor on $\F_p^n$ of complexity $(\ell',q')$ and rank at least $\rho(\ell+q)$, satisfying $\calL'\supseteq \calL$, $\calQ'\supseteq \calQ$, and $\ell'+q'\leq \ell+q+\frac{1}{3}\log_p(\e^{-1/32})$.  Let $(u,v)\in \F_p^{\ell}\times \F_p^{q}$ be such that $\Sigma_{\calB}(d)=(u,v)$. Then for any $(u',v')\in \F_p^{\ell'-\ell}\times \F_p^{q'-q}$, 
 $$
 \Big|\frac{|A\cap B(u;v)|}{|B(u;v)|}-\frac{|A\cap  B'(uu';vv')|}{|B'(uu';vv')|}\Big|\leq 5\e^{1/32}.
 $$ 
\end{proposition}
 
We now prove the main result of this section, which clearly implies Theorem \ref{thm:lowerbound}. 

\begin{theorem}
For all primes $p>2$, there exist $\e\in (0,1)$ and a polynomial growth function $\rho$ such that for all sufficiently large odd integers $n$,  there exists $A\subseteq\F_p^n$ satisfying the following.

Fix integers $\ell',q'\geq 1$, and suppose $\calB'=(\calL',\calQ')$ is a quadratic factor on $\F_p^n$ of complexity $(\ell',q')$ and rank at least $\rho(\ell'+q')$ such that there exists $\Gamma\subseteq \F_p^{3\ell'+6q'}$ satisfying 
\begin{enumerate}[label=\normalfont(\roman*)]
\item $|\Gamma|\geq (1-\e)p^{3\ell'+6q'}$;
\item for all $d\in \Gamma$, $\|1_A-\alpha_{B(\Sigma(d))}\|_{U^3(d)}<\e$.
\end{enumerate}
Then $q'\geq (1-\e)\e^{-1/64}$.
\end{theorem}
\begin{proof}
Fix $0<\e<10^{-32}$ sufficiently small.  Let $\rho_0=\rho_0(p,\e^4)$ from Corollary \ref{cor:sizeofatoms}, and let $\rho_1=\rho_1(p,\e)$ be from Proposition \ref{prop:preservedensity}.  Now let $\rho$ be the growth function defined by setting 
$$
\rho(x)=\max\{ \rho_0(x+\lceil \log_p(\e^{-1/64})\rceil+1),\rho_1(x+\lceil\log_p(\e^{-1/64})\rceil+1)\}.
$$
We note $\rho$ is bounded above by a polynomial in $x$ and $\e^{-1}$ because $\rho_0$ and $\rho_1$ are.  Let $n$ be a sufficiently large odd integer.   

We begin by defining our subset $A\subseteq \F_p^n$.  Let $\ell=\lceil \log_p(\e^{-1/64})\rceil$ and $q=p^{\ell}$.  Let $\{M_1,\ldots, M_{p^{\ell}}\}$ be a purely quadratic factor on $\F_p^n$ of rank $n$ (this exists by \cite[Lemma 3]{Dumas.2011}).   Define $\calL=\{e_1,\ldots, e_{\ell}\}$, where $e_1,\ldots, e_{\ell}$ are the first $\ell$ basis vectors in $\F_p^n$, and let $\calQ=\{N_c: c\in \F_p^{\ell}\}$ be any re-indexing of $\{M_1,\ldots, M_{p^{\ell}}\}$ by elements of $\F_p^{\ell}$.  We define our set $A$ as follows.
$$
A=\bigcup_{c\in \F_p^{\ell}}\Big(L(c)\cap  \{x: x^TN_cx=0\}\Big).
$$
Suppose now $\ell',q'\geq 1$ are integers, $\calB'=(\calL',\calQ')$ is a quadratic factor on $\F_p^n$ of complexity $(\ell',q')$, and $\Gamma\subseteq \F_p^{3\ell'+6q'}$ with the properties
\begin{enumerate}
\item $\calB'$ has rank at least $\rho(\ell'+q')$;
\item $|\Gamma|\geq (1-\e)p^{3\ell'+6q'}$;
 \item for all $d\in \Gamma$, $\|1_A-\alpha_{B(\Sigma(d))}\|_{U^3(d)}<\e$.  
 \end{enumerate}
Note that the rank of $\calB'$ is at most $n$, so $n\geq \rho(\ell'+q')$. Let
$$
\Sigma(\Gamma)=\{b\in \F_p^{\ell'+q'}: \text{ for some }d\in \Gamma, b=\Sigma_{\calB}(d)\}.
$$
For any $d\in \Gamma$, there are at $p^{2\ell'+5q'}$ many $d'\in \Gamma$ with $\Sigma_{\calB'}(d)=\Sigma_{\calB'}(d')$, so 
$$
|\Sigma(\Gamma)|\geq p^{-2\ell'-5q'}|\Gamma|\geq (1-\e)p^{\ell'+q'},
$$
where the second inequality uses that $|\Gamma|\geq (1-\e)p^{3\ell'+6q'}$. Let 
$$
X=\bigcup_{d\in \Gamma}B'(\Sigma(d))=\bigcup_{b\in \Sigma(\Gamma)}B'(b).
$$
By Corollary \ref{cor:sizeofatoms} and the above lower bound on $\Sigma(\Gamma)$, we have   
$$
|X|\geq |\Sigma(\Gamma)|(1-\e^4)p^{n-\ell'-q'}\geq  (1-\e^4)(1-\e)p^n\geq (1-2\e)p^{n}.
$$
Define
$$
\Theta=\{c\in \F_p^{\ell}: X\cap L(c)\neq \emptyset\}.  
$$
Note that 
$$
(1-\e)p^n\leq |X|\leq \sum_{c\in \Theta}|X\cap L(c)|\leq  \sum_{c\in \Theta}|L(c)|=|\Theta| p^{n-\ell}.
$$
Consequently, $|\Theta|\geq (1-\e)p^{\ell}$.

\begin{claim}\label{cl:c} For all $c\in \Theta$, $\calQ'\cup \{N_c\}$ has rank less than $\rho(\ell'+q')$.
\end{claim}
\begin{proof}
Fix $c\in \Theta$, and suppose towards a contradiction $\calQ'\cup \{N_c\}$ has rank at least $\rho(\ell'+q')$.  By definition of $c\in \Theta$, $L(c)\cap X\neq \emptyset$, meaning there exists some $b\in \Sigma(\Gamma)$ such that $L(c)\cap B'(b)\neq \emptyset$.   

Let $\calL''$ be a minimal linear factor on $\F_p^n$ satisfying $\calL'\subseteq \calL''$ and $\calL\cup \calL'\subseteq \Span(\calL'')$.  Now set $\calQ''=\calQ'\cup \{N_c\}$, and consider the quadratic factor $\calB''=(\calL'',\calQ'')$.  Note $\calB''$ has complexity $(\ell'',q'')$ where $\ell'':=|\calL''|\leq \ell'+\ell$ and $q''=q'+1$.  Therefore,
$$
\ell''+q''\leq \ell'+\ell+q'+1\leq \ell'+q'+\lceil\log_p(\e^{-1/64})\rceil+1\leq \ell'+q'+\frac{1}{3}\log_p(\e^{-1/32}),
$$    
where the second inequality is by our assumption on $\ell$, and the last is because $\e$ is sufficiently small.  By our assumption on the rank of $\calB'$ and our choice of $\rho$, we can conclude 
$$
\rk(\calB'')=\rk(\calB')\geq \rho(\ell'+q') \geq \max\{\rho_0(\ell''+q''),\rho_1(\ell''+q'')\}.
$$
Observe that since since $L(c)\cap B'(b)\neq \emptyset$, it must be an atom of $(\calL'',\calQ')$.  Consequently, we have that for any $u\in \F_p$, the set 
$$
L(c)\cap B'(b)\cap \{x\in \F_p^n: x^TN_cx=u\}
$$
 is an atom of $\calB''$ (note it is non-empty by by Corollary \ref{cor:sizeofatoms}). By Proposition \ref{prop:preservedensity},  the following hold.  
$$
\Big| \frac{|A\cap L(c)\cap B(d)\cap \{x: x^TN_cx=0\}|}{|B(d)\cap L(c)\cap \{x: x^TN_cx=0\}|}-\frac{|A\cap B'(d)|}{|B'(d)|}\Big|\leq 5\e^{1/32}
$$
and
$$
\Big| \frac{|A\cap L(c)\cap B(d)\cap \{x: x^TN_cx=1\}|}{|B(d)\cap L(c)\cap \{x: x^TN_cx=1\}|}-\frac{|A\cap B'(d)|}{|B'(d)|}\Big|\leq 5\e^{1/32},
$$
and thus
$$
\left| \frac{|A\cap B(d)\cap L(c)\cap \{x: x^TN_cx=0\}|}{|B(d)\cap L(c)\cap \{x: x^TN_cx=0\}|}- \frac{|A\cap B(d)\cap L(c)\cap \{x: x^TN_cx=1\}|}{|B(d)\cap L(c)\cap \{x: x^TN_cx=1\}|}\right|\leq 10\e_1^{1/32}.
$$
However, by definition of $A$, one of these values is $1$ while the other is $0$, and thus the above inequality is impossible, as $10\e^{1/32}<1$, so we have arrived at a contradiction.  This finishes the proof of the claim.
\end{proof}

The above claim tells us that for all $c\in \Theta$, $\calQ'\cup \{N_c\}$ has rank less than $\rho(\ell'+q')$, and consequently, there are some $\rho_c,\lambda_1^c,\ldots, \lambda_{q'}^c\in \F_p$, not all $0$, so that 
$$
\rk(\rho_c N_c+\sum_{i=1}^{q'}\lambda_i^c Q_i)<\rho(\ell'+q').
$$
We observe that since $\calQ'$ has rank at least $\rho(\ell'+q')$, we must have $\rho_c\neq 0$.  Further, since $N_c$ has rank $n$, there is also some $i\in [q']$ such that $\lambda_i^c\neq 0$.   

Let $\mathbf{V}$ be the $\F_p$-vector space consisting of all $n\times n$ matrices with entries in $\F_p$. Consider the following subset of $\mathbf{V}$: 
$$
\calS=\Big\{\sum_{i=1}^{q'}\lambda_i^cQ_i: c\in \Theta\Big\}.
$$
We claim that $\calS$  is linearly independent as a subset of $\mathbf{V}$.  Indeed, suppose this was false.  Then there exists a tuple $(\mu_c)_{c\in \Theta}$ of elements of $\F_p$, not all entries of which are $0$, such that  
$$
\sum_{c\in \Theta}\mu_c\Big(\sum_{i=1}^{q'}\lambda_i^cQ_i\Big)=0
$$
This implies
\begin{align}\label{theta}
  \sum_{c\in \Theta}\mu_c( \rho_cN_c+\sum_{i=1}^{q'}\lambda_i^cQ_i)=\sum_{c\in \Theta}\mu_c\rho_cN_c.
\end{align}
However, by assumption, the left hand side of (\ref{theta}) is a sum of matrices, each of rank less than $\rho(\ell'+q')$. Thus, by subadditivity, the left hand side of (\ref{theta}) has rank less than $\rho(\ell'+q')$. However, each $\rho_c$ is nonzero, and some $\mu_c$ is nonzero, meaning some coefficient $\mu_c\rho_c$ on the right hand side of (\ref{theta}) is nonzero.  Since $\{N_c: c\in \F_p^{\ell}\}$ has rank $n$, this implies the right hand side of (\ref{theta}) has rank $n$.  We have now arrived at a contradiction, since $n\geq \rho(\ell'+q')$.   This finishes our verification that $\calS$ is a linearly independent set in $\mathbf{V}$.

Because $\calS$ is linearly independent in $\mathbf{V}$,  the following holds, where the subspaces are computed in $\mathbf{V}$:
\begin{align}\label{theta1}
p^{|\Theta|}=p^{|\calS|}=|\langle \calS\rangle|\leq |\langle \calQ'\rangle|\leq p^{q'},
\end{align}
where the first inequality is because $\calS\subseteq \langle \calQ'\rangle$ by definition, and the last inequality is because $|\calQ'|=q'$.  From this we can conclude $q'\geq |\Theta|$. On the other hand, the lower bound for $\Theta$ obtained before Claim \ref{cl:c} tells us $|\Theta|\geq (1-\e)p^{\ell}$. Thus $q'\geq (1-\e)p^{\ell}\geq (1-\e)\e^{-1/64}$, which finishes the proof.
\end{proof}


\providecommand{\bysame}{\leavevmode\hbox to3em{\hrulefill}\thinspace}
\providecommand{\MR}{\relax\ifhmode\unskip\space\fi MR }
\providecommand{\MRhref}[2]{%
  \href{http://www.ams.org/mathscinet-getitem?mr=#1}{#2}
}
\providecommand{\href}[2]{#2}

\end{document}